\documentclass[12pt,preprint]{elsarticle}

\usepackage{amsmath,amssymb}

\usepackage[utf8]{inputenc}
\usepackage[english]{babel}

\usepackage{graphicx}
\usepackage{xcolor}
\usepackage{enumitem}
\usepackage{float}
\usepackage{subcaption}

\usepackage{hyperref}

\biboptions{sort&compress}

\newcounter{rem}
\newtheorem{remark}[rem]{Remark}

\def\dtp#1{\mathop {#1}\limits_{+\tau}}
\def\dtm#1{\mathop {#1}\limits_{-\tau}}
\def\dsp#1{\mathop {#1}\limits_{+s}}
\def\dsm#1{\mathop {#1}\limits_{-s}}

\makeatletter
\def\ps@pprintTitle{%
  \let\@oddhead\@empty
  \let\@evenhead\@empty
  \def\@oddfoot{\rightline{\emph{\today}}}
  \let\@evenfoot\@oddfoot
}
\makeatother

\begin{document}

\begin{frontmatter}

\title{Conservative invariant finite-difference schemes for the
 modified shallow water equations in Lagrangian coordinates}

\author[Keldyshaddress]{V.~A. Dorodnitsyn}
\ead{Dorodnitsyn@keldysh.ru,dorod2007@gmail.com}

\author[SUTaddress]{E.~I. Kaptsov}
\ead{evgkaptsov@gmail.com}

\author[SUTaddress]{S.~V. Meleshko\corref{correspondingauthor}}
\cortext[correspondingauthor]{Corresponding author}
\ead{sergey@math.sut.ac.th}

\address[Keldyshaddress]{Keldysh Institute of Applied Mathematics,\\
Russian Academy of Science, Miusskaya Pl. 4, Moscow, 125047, Russia}

\address[SUTaddress]{School of Mathematics, Institute of Science, \\
Suranaree University of Technology, Nakhon Ratchasima, 30000, Thailand}

\begin{abstract}
 The one-dimensional modified shallow water equations
in Lagrangian coordinates are considered. It is shown the
relationship between symmetries and conservation laws in Lagrangian coordinates,
in mass Lagrangian variables, and Eulerian coordinates.
For equations in Lagrangian coordinates an invariant finite-difference scheme is constructed
for all cases for which conservation laws exist in the differential model.
Such schemes possess the difference analogues of the conservation laws of mass, momentum,
energy, the law of center of mass motion for horizontal, inclined and parabolic bottom topographies.
Invariant conservative difference scheme is tested numerically in comparison with naive
approximation invariant scheme.
\end{abstract}


\end{frontmatter}

\bigskip

%\cite{DORODNITSYN2019201,KOZLOV2019,bk:Korobitsyn_scheme[1989]}

\section{Introduction}

Mathematical modeling of physical phenomena is one of the main
streams in continuum mechanics. Such phenomena  as hydraulic
currents, coastal currents, currents in rivers and lakes, currents
in water intakes, technical troughs and trays, tsunami simulation,
propagation breakthrough of waves and tidal pine forests in rivers,
the spread of heavy gases and impurities in the atmospheres of
planets, atmospheric movements scales used in weather forecasting
require mathematical consideration.

Motion of ideal fluid flow under the force of gravity can be modeled
by means of the Euler equations. However, the full Euler equations,
even under the assumption of incompressibility, barotropy and
absence of rotation, are still rather complicated for describing
waves on a surface. One of these difficulties is that  the free
surface is a part of the solution. This difficulty has motivated
scientists to derive simpler equations. For this reason development
of approximate models and their analysis by analytical and numerical
methods is an actual  problem.

The need to reduce the original equations to simpler equations led to the construction of
asymptotic expansion models with respect to a small parameter determined by
the ratio of the depth of the fluid to the characteristic linear size.
One class of such equations is the class of shallow water equations.
There are many approaches for deriving shallow water models, a review
of which can be found in \cite{bk:BonnetoBarthelemy_et[2011],bk:KhakimzyanovDutykhFedotovaMitsotakis}.

The shallow water equations describe the motion of incompressible fluid
in the gravitational field if the depth of the liquid layer is small enough.
They are widely used in the description of processes in the atmosphere,
water basins, modeling of tidal oscillations, tsunami waves and gravitational
waves~(see the classical papers such as~\cite{bk:Whitham[1974],bk:Ovsyannikov[2003]} and
detailed description in, for example,~\cite{bk:PetrosyanBook[2010],bk:Vallis[2006]}).

For solving real-world problems it is also necessary to consider additional impacts determined by specific flow conditions. This leads to the appearance of additional terms in the system of shallow water equations. In particular, shallow water flows on plane surface in the presence of the weak vertical inhomogeneities in the
initial conditions contain additional terms appearing as the result of depth averaging of the nonlinear terms in the initial fluid equations \cite{bk:LeVeque,bk:Pedlosky}.
 One of the models including this effect in consideration is proposed in \cite{bk:KarelskyPetrosyan2006,bk:KarelskyPetrosyan2008}:
\begin{subequations} \label{modSW_Eulerian}
\begin{equation}
\frac{\partial \rho}{\partial t} + \frac{\partial}{\partial x}(\rho u) = 0,
\end{equation}
\begin{equation}
\frac{\partial u}{\partial t} + u \frac{\partial u}{\partial x}
    + g \left(1 + \frac{g_1}{\rho}\right) \frac{\partial \rho}{\partial x} - g H^\prime(x) = 0,
\end{equation}
\end{subequations}
where $t$ is time, $x$ is the Eulerian space coordinate, $u$ is the
velocity, $\rho$ is the depth of the layer of fluid,
the differentiable function~$H(x)$ describes the bottom topography, $g$ is the
gravitational acceleration, and~$g_1/\rho$ is the additional term. The
latter term  describes an advective transport of impulse as a result
of the dependence of horizontal shallow water flows on vertical
coordinate \cite{bk:KarelskyPetrosyan2006}. It is assumed that $g_1
\neq 0$ in this paper,  otherwise equations~(\ref{modSW_Eulerian})
become standard one-dimensional shallow water equations.
We use the letter~$\rho$ for the depth of the fluid instead of usually used~$h$ to differ it from the standard notation for the finite-difference mesh spacing.

It is wellknown that symmetry of mathematical model is intrinsic
property inherited from physical phenomena.    One of the tools for
studying symmetries is Lie group analysis
\cite{bk:Ovsyannikov[1982],bk:Olver[1986]}  which is a basic method
for constructing exact solutions of ODEs and partial differential
equations. Even in case of the one-dimensional shallow water
equations for the flat bottom one meets certain difficulties to
obtain nontrivial exact solutions. Applications of Lie groups  to
differential equations
 is the subject of many books and review
articles~\cite{bk:Ovsyannikov[1982],bk:Olver[1986],bk:Ibragimov1985,bk:Bluman1989,
bk:HandbookLie_v1,bk:Gaeta1994}.

The group properties of the shallow water equations were studied
in numerous papers (see~\cite{bk:HandbookLie_v2,bk:LeviNicciRogersWint[1989],bk:ClarksonBila[2006],bk:Andronikos2019}).
Group classification and first integrals of these equations can
be found in~\cite{bk:AksenovDruzkov_classif[2019],bk:KaptsovMeleshko_1D_classf[2018]}.
It was shown~(see, e.~g.,~\cite{bk:SiriwatKaewmaneeMeleshko2016,bk:KaptsovMeleshko_1D_classf[2018]}) that the shallow
water equations in Lagrangian coordinates can be obtained
as Euler--Lagrange equations of Lagrangian functions of a special kind.
Some exact solutions also can be found in~\cite{bk:PetrosyanBook[2010],bk:Bernetti[2008],bk:HanHantke[2012]}.

More recently applications of Lie groups have been extended
to difference equations~\cite{Maeda1, Maeda2,Dor_1, Dor_2, Dor_3,
bk:DorodKozlovWint[2004],KimOlv04,KimOlv07,[LW-2],bk:DorodKozlovWinternitz[2000],bk:Dorod_NoetherTh[2001],
Quisp, [LW-3], bk:Dorodnitsyn[2011], Vinet,
bk:Hydon_book[2014],bk:DorodKozlovWintKaptsov[2015]}.
%%%%%%%%%%%%%%%%%%%%%%%%%%%%%%%%%%%%%%%%%%%

The applications have certain peculiarities related with nonlocal
character of difference operators and geometrical structure of
difference mesh~\cite{bk:Dorodnitsyn[2011]}. The results were
obtained both for the ordinary and partial difference equations and
systems~(e.~g.,~\cite{Dor_1,Dor_2,
[Dheat],bk:Dorodnitsyn[2011],bk:DorodKaptsov[2013],bk:DorodKaptsov_Ermakov[2016]}).
The finite-difference analogues of
Lagrangian~\cite{Dor_3,bk:DorodKozlovWint[2004],bk:Dorod_Hamilt[2011]}
and the
Hamiltonian~\cite{bk:Dorod_Hamilt[2010],bk:Dorod_Hamilt[2011]}
formalism were developed. Moreover, it was shown that in case of an
absence of the Lagrangian and Hamiltonian there is a more general
approach based on the Lagrange operators identity and adjoint
equations method~\cite{bk:BlumanAnco2002}, the difference analogue of which
was developed in
\cite{bk:DorodKozlovWintKaptsov[2014],bk:DorodKozlovWintKaptsov[2015]}.
In all three approaches the starting point is the symmetry of the
differential equations, preserved in difference equations and
meshes. In this article we follow Lagrangian approach and the
Noether theorem. We also discuss so called direct
method~\cite{bk:Bluman1997,bk:BlumanAnco2002,bk:BlumanCheviakovAnco} for
constructing conservation laws.

The invariant finite-difference schemes   are of  particular
interest as far it should preserve sym\-met\-ries and since the
geometric properties of the original equations, in particular it
have symmetry reductions on subgroups and exact invariant solutions.
In this paper we mostly concentrate on invariant schemes in the
Lagrange coordinate system.

The present paper is devoted to the construction of invariant
conservative difference schemes for the modified shallow water
equations (\ref{modSW_Eulerian}) in Lagrangian coordinates and mass
Lagrangian coordinates. The base of the construction is the
invariant scheme  already developed for standard shallow water model~\cite{dorodnitsyn2019shallow}.

\bigskip

The paper is organized as follows.
In Section~\ref{sec:diffeqns}, the modified shallow water equations in Lagrange coordinates are given. It is shown that they can be obtained as the Euler--Lagrange equations for a certain Lagrangian.
Symmetries and conservation laws for the modified shallow water equations for various bottom topographies
are given in Lagrangian, mass Lagrangian and Eulerian coordinates in Section~\ref{sec:symsAndCLs}.
In Section~\ref{sec:discr}, the equations in Lagrangian coordinates are discretized. The constructed finite-difference schemes for various bottom topographies are invariant and possess finite-difference analogues of the differential conservation laws.
In addition, a scheme defined on a reduced finite-difference stencil is constructed in mass Lagrangian coordinates.
Approaches to the numerical comparison of various schemes in Lagrangian coordinates are also discussed in the section.
In Section~\ref{sec:numimpl}, the numerical implementation of the constructed schemes is carried out.
The obtained schemes and numerical results are discussed in Conclusion.

\section{One-dimensional modified shallow water equations in Lagrangian coordinates}
\label{sec:diffeqns}

Modelling physical phenomena in continuum mechanics is considered
in two distinct ways. The typical approach uses Eulerian coordinates,
where flow quantities at each instant of time during motion are described
at fixed points. Alternatively, the Lagrangian description is used,
where the particles are identified by the positions which they occupy
at some initial time.

Following the classical gas dynamics equations, introduce Lagrangian
variables $(\xi,t)$, relating the Eulerian variables $(x,t)$ and
Lagrangian variables by the equation
\[
x=\tilde{\varphi}(t,\xi),
\]
where the function $\tilde{\varphi}(t,\xi)$ satisfies the Cauchy
problem
\[
\tilde{\varphi}_{t}=u(t,\tilde{\varphi}),\,\,\,\tilde{\varphi}(t_0,\xi)=\xi.
\]
Subindex of a function $f$ denotes a corresponding partial derivative.

Let the dependent variables in Lagrangian variables are denoted as
$\tilde{\rho}(t,\xi)$ and $\tilde{u}(t,\xi)$. Their relations with
counterparts are given by the formulas
\[
\tilde{\rho}(t,\xi)=\rho(t,\tilde{\varphi}(t,\xi)),\,\,\,\tilde{u}(t,\xi)=u(t,\tilde{\varphi}(t,\xi)).
\]
Differentiating the latter with respect to $\xi$ and $t$, one finds
\[
\rho_{x}=\tilde{\rho}_{\xi}/\tilde{\varphi}_{\xi},\,\,\,\rho_{t}+u\rho_{x}=\tilde{\rho}_{t},\,\,\,
u_{x}=\tilde{u}_{\xi}/\tilde{\varphi}_{\xi},\,\,\,u_{t}+uu_{x}=\tilde{u}_{t}.
\]
Noticing that $\tilde{u}_{\xi}=\tilde{\varphi}_{t\xi}$, the continuity
equation becomes
\[
(\tilde{\rho}\tilde{\varphi}_{\xi})_{t}=0.
\]
Hence,
\[
\tilde{\rho}(t,\xi)\tilde{\varphi}_{\xi}(t,\xi)=\tilde{\rho}_{0}(\xi),
\]
where $\tilde{\rho}_{0}(\xi)=\tilde{\rho}(\xi,0)$. Introducing the change
$\xi=\alpha(s)$, where
\[
\alpha^{\prime}(s)\tilde{\rho}_{0}(\alpha(s))=1,
\]
one derives the independent variables $(s,t)$, which are called in
the classical gas dynamics by the mass Lagrangian coordinates. In
the mass Lagrangian coordinates one obtains that
\[
\bar{\rho}(t,s)=\frac{1}{\bar{\varphi}_{s}(t,s)},
\]
where
\[
\bar{\rho}(t,s)=\tilde{\rho}(t,\alpha(s)),\,\,\,\bar{u}(t,s)=\tilde{u}(t,\alpha(s)),\,\,\,
\bar{\varphi}(t,s)=\tilde{\varphi}(t,\alpha(s)).
\]
Further the sign bar `~$\bar{\ensuremath{}}$ ' is omitted.

In the mass Lagrangian coordinates, the second equation of system~(\ref{modSW_Eulerian})
becomes
\begin{equation}\label{modSW_EulerianOrig}
\varphi_{tt}-g\frac{\varphi_{ss}}{\varphi_{s}^{3}}-gg_{1}\frac{\varphi_{ss}}{\varphi_{s}^{2}}
- g H^\prime(\varphi)=0.
\end{equation}

One of the advantages of using the mass Lagrangian coordinates is
the existence of a Lagrangian whose Euler-Lagrange equations is equivalent
to equation (\ref{modSW_EulerianOrig}). In Lagrangian mechanics the
problem of determining whether a given system of differential equations
can arise as the Euler-Lagrange equations for some Lagrangian function
is called the Helmholtz problem.

For finding a Lagrangian for which equation (\ref{modSW_EulerianOrig})
is the Euler-Lagrange equations one has to solve the following problem.
Let ${\cal L}(t,s,\varphi,\varphi_{s},\varphi_{t})$ be a corresponding
Lagrangian. Then, substituting ${\cal L}$ into the equation
\begin{equation}
{\displaystyle \frac{\delta{\cal L}}{\delta\varphi}}=0,\label{eq:oct08.5}
\end{equation}
excluding the derivative $\varphi_{tt}$ found from equation (\ref{modSW_EulerianOrig}),
and splitting them with respect to the parametric derivatives $\varphi_{ts},\,\,\,\varphi_{ss}$,
one obtains an overdetermined system of equations for the function
${\cal L}$. Here ${\displaystyle \frac{\delta}{\delta\varphi}}$
is the variational derivative:
\[
{\displaystyle \frac{\delta F}{\delta\varphi}=\frac{\partial F}{\partial\varphi}-D_{t}^{L}\left(\frac{\partial F}{\partial\varphi_{t}}\right)-D_{s}\left(\frac{\partial F}{\partial\varphi_{s}}\right),}
\]
where $F$ is an arbitrary function, $D_{t}^{L}$, $D_{s}$ are total
derivatives with respect to $t$ and $s$, respectively.
It supposed that
\begin{equation}
\Delta=\frac{\partial^{2}{\cal L}}{\partial\varphi_{t}^{2}}\neq0.\label{eq:oct08.6}
\end{equation}
Any such solution of this overdetermined system of equations satisfying
condition (\ref{eq:oct08.6}) gives a sought Lagrangian.

The calculations give that
\[
\mathcal{L}=k\left(\frac{\varphi_{t}^{2}}{2}+gg_{1}\ln(\varphi_{s})-g\varphi_{s}^{-1} + g H(\varphi)\right)
 + G,
\]
where $G=\varphi_{t}L_{1}+\varphi_{s}L_{2}+L_{3}$, and the functions
$L_{i}(t,s,\varphi),\,\,\,(i=1,2,3)$ satisfy the condition
\[
L_{1t}+L_{2s}=L_{1\varphi}+L_{2\varphi}.
\]
 Noticing that ${\displaystyle \frac{\delta G}{\delta\varphi}}=0, $and
because of the condition (\ref{eq:oct08.6}), one derives that the
seeking Lagrangian is can be chosen as
\begin{equation}
\label{eq:nov12.1}
{\cal L}=\frac{\varphi_{t}^{2}}{2}+gg_{1}\ln(\varphi_{s})-g\varphi_{s}^{-1} + g H(\varphi).
\end{equation}

By means of the transformation
\[
t = \tilde{t}/\sqrt{2},
\qquad
s= \tilde{s}/g,
\qquad
g_{1}=2\gamma_{1}/g,
\qquad
H = 2(\tilde{H} - \gamma_1 \ln g)/g,
\]
the Lagrangian (\ref{eq:nov12.1}) and equation~(\ref{modSW_EulerianOrig})
are brought to the forms\footnote{Here and further the symbol~$\tilde{\,}$ is also omitted.}
\begin{equation}
\label{eq:nov12.1-1}
{\cal L}=\frac{\varphi_{t}^{2}}{2}+\gamma_{1}\ln(\varphi_{s})-\frac{1}{2\varphi_{s}} + H(\varphi),
\end{equation}
\begin{equation}\label{modSW_Lagr}
\varphi_{tt}-\frac{\varphi_{ss}}{\varphi_{s}^{3}}-\gamma_{1}\frac{\varphi_{ss}}{\varphi_{s}^{2}}
- H^\prime(\varphi)
=0.
\end{equation}

%%%%%%%%%%%%%%%%%%%%%%%%%%%%%%%%%%%%%%%%%%%%%%%%%%%%%
%%%%%%%%%%%%%%%%%%%%%%%%%%%%%%%%%%%%%%%%%%%%%%%%%%%%%
% TO BE MODIFIED BY E.K.
%%%%%%%%%%%%%%%%%%%%%%%%%%%%%%%%%%%%%%%%%%%%%%%%%%%%%
%%%%%%%%%%%%%%%%%%%%%%%%%%%%%%%%%%%%%%%%%%%%%%%%%%%%%

\section{Symmetries and conservation laws of the modified shallow water equations}
\label{sec:symsAndCLs}

\subsection{Symmetries and conservation laws in Lagrangian coordinates}

In~\cite{bk:KaptsovMeleshko_1D_classf[2018]} the group classification
of the one-dimensional Euler--Lagrange equations of continuum mechanics was carried out by the authors.
In particular, it was shown that the modified shallow water equation~(\ref{modSW_Lagr})
is the Euler--Lagrange equation for the Lagrangian
\begin{equation} \label{lagr}
\mathcal{L} = \frac{\varphi_t^2}{2} - \frac{1}{2\varphi_s} + \gamma_1 \ln \varphi_s + H(\varphi).
\end{equation}

%\todo{Discuss notation: $H^\prime(\varphi)=0$ vs $H(\varphi)=\text{const}$ vs $H(x)=\text{const}$.}

The group classification states that in case~$H^\prime(\varphi)$ is arbitrary,
equation~(\ref{modSW_Lagr}) admits two generators
\begin{equation} \label{LagrXkern}
X_1 = \frac{\partial}{\partial t},
\qquad
X_2 = \frac{\partial}{\partial s}.
\end{equation}
In case $H^\prime(\varphi) = 0$, there is the following extension of the admitted Lie algebra
\begin{equation} \label{SymExt1}
X_3 = \frac{\partial}{\partial \varphi},
\qquad
X_4 = t\frac{\partial}{\partial \varphi},
\qquad
X_5 = t\frac{\partial}{\partial t}
+ s\frac{\partial}{\partial s}
+ \varphi \frac{\partial}{\partial \varphi}.
\end{equation}

It is shown below that the case of an inclined bottom~$H^\prime(\varphi)=\text{const}$ is reduced to the case of a horizontal bottom by a point transformation. Therefore, we do not consider this case here.

In case of a parabolic bottom~$H^\prime(\varphi) = \varphi$,
the extension of the admitted Lie algebra is~\cite{bk:DorKapJMPSW2021}
\begin{equation} \label{SymExt2}
X^+_3 = e^t \frac{\partial}{\partial \varphi},
\qquad
X^+_4 = e^{-t} \frac{\partial}{\partial \varphi}.
\end{equation}
In case $H^\prime(\varphi) = -\varphi$, one has
\begin{equation} \label{SymExt3}
X^-_3 = \sin t\, \frac{\partial}{\partial \varphi},
\qquad
X^-_4 = \cos t \, \frac{\partial}{\partial \varphi}.
\end{equation}
In case $H^\prime = 1/\varphi$,  the extension only consists of the generator
\begin{equation}
X_5 = t\frac{\partial}{\partial t}
+ s\frac{\partial}{\partial s}
+ \varphi \frac{\partial}{\partial \varphi}.
\end{equation}

\par
\bigskip
\par

Usually the Lagrangian admits some subset of generators of the Lie algebra admitted by the corresponding Euler--Lagrange equations. Such generators are called variational symmetries. The presence of variational symmetries greatly simplifies the search for conservation laws for the Euler--Lagrange equations.
This can be done using Noether's theorem~\cite{bk:Noether1918,bk:Ibragimov1985} which establishes a connection between variational symmetries and conservation laws.

Recall that (local) conservation laws of the system
\begin{equation}\label{genPDEsys}
F^i(t,s,\varphi, \varphi_t, \varphi_s, \varphi_{tt}, \varphi_{ts}, \varphi_{ss})=0,
\qquad
i = 1, 2, ..., n,
\end{equation}
can be represented in the form
\begin{equation}
(T^t)_t + (T^s)_s = \Lambda_j F^j = 0,
\end{equation}
% іќііі-і-і-ііі- TЂі-TЃTЃііі-іїі-T‚TЊ іїTЂі- іїTЂTЏі-і-іїіі і-ііT‚і-і+ T‚TѓT‚
where
$T^t=T^t(t,s,\varphi, \varphi_t, \varphi_s)$ is a conserved density,
$T^s=T^s(t,s,\varphi, \varphi_t, \varphi_s)$ is a conserved flux,
and~$\Lambda_j = \Lambda_j(t,s, \varphi, \varphi_t, \varphi_s)$, $j=1,...,n$
are so called conservation law multipliers~\cite{bk:BlumanAnco2002}.
Further on~$\Lambda_{\alpha j}$ denote the multipliers corresponding
to the conservation law obtained using Noether's theorem for the generator~$X_\alpha$.
For brevity, in case $n = 1$, the subscript~$j$ is omitted.

\smallskip

There is a one-to-one correspondence between the equivalence classes of conservation laws and conservation law multipliers~\cite{bk:Olver[1986]}, and it is often easier to find the multipliers first.
Conservation law multipliers for a given system~(\ref{genPDEsys}) can be found by means of direct method~\cite{bk:Bluman1997,bk:BlumanAnco2002}.
The method consists of applying the Euler operator
\begin{equation}
\displaystyle
E_\varphi \equiv
    \frac{\partial}{\partial \varphi}
        - D_i \frac{\partial}{\partial \varphi_i}
        + \cdots
        + (-1)^s D_{i_1} \cdots D_{i_s} \frac{\partial}{\partial \varphi_{i_1 \cdots i_s}} + \cdots
\end{equation}
to the expression $\Lambda_j F^j$ for unknown multipliers.
(Here $D_i$ denotes the total differentiation operator with respect to~$i$-th independent variable.)
Since the Euler operator identically vanishes divergence expressions,
to find conservation law multipliers one has to solve the equation
\begin{equation}
E_\varphi(\Lambda_j F^j) \equiv 0.
\end{equation}

\smallskip

In the present section, Noether's theorem is used. It turned out that
the finite-difference analogue of the direct method is effective when constructing
conservative finite-difference schemes. This is discussed in the related sections below.

\bigskip

By means of Noether's theorem one finds the following conservation laws,
and the corresponding conservation law multipliers of equation~(\ref{modSW_Lagr}) by the direct method.

\smallskip

For an arbitrary differentiable function~$H(x)$ and the generators~$X_1$ and~$X_2$ one obtains
the conservation law of energy
\begin{equation}\label{ws_cl_01}
\def\arraystretch{1.75}
\Lambda_1 = \varphi_t,
\qquad
\left(
    \frac{\varphi_t^2}{2}
    + \frac{1}{2 \varphi_s}
    - \gamma_1 \ln \varphi_s
    - H
\right)_t
+ \left(
       \varphi_t \left(
            \frac{\gamma_1}{\varphi_s}
            + \frac{1}{2 \varphi_s^2}
    \right)
\right)_s = 0,
\end{equation}
and the conservation law of momentum
\begin{equation}\label{ws_cl_02}
\Lambda_2 = \varphi_s,
\qquad
\left( \varphi_t \varphi_s \right)_t
    + \left(
        \frac{1}{\varphi_s}
        - \frac{1}{2} \varphi_t^2
        -\gamma_1 \ln \varphi_s
        - H
        \right)_s = 0.
\end{equation}

\smallskip

Conservation laws obtained for $H^\prime=0$ and the generators $X_3$ and $X_4$
are the alternative form of the momentum conservation law
\begin{equation}\label{ws_cl_03}
 \Lambda_3 = 1,
  \qquad
        (\varphi_t)_t
            + \left(
                \frac{1}{2 \varphi_s^2}
                + \frac{\gamma_1}{\varphi_s}
            \right)_s = 0,
\end{equation}
and the center of mass law
\begin{equation}\label{ws_cl_04}
        \Lambda_4 = t,
    \qquad
        \left( t \varphi_t - \varphi \right)_t
            + \left(
                \frac{t}{2 \varphi_s^2}
                + \frac{\gamma_1 t}{\varphi_s}
        \right)_s = 0.
\end{equation}
%The generator $X^1$ does not produce a conservation law since it does not satisfy Noether's theorem.

\smallskip

The conservation laws corresponding to $H^\prime = \varphi$ and the generators $X^{+}_3$ and $X^{+}_4$ are
\begin{equation} \label{ws_cl_05-1}
\Lambda^{+}_3 = e^t,
\qquad
\left( e^t ({\varphi} - {\varphi}_t)\right)_t
    - \left(e^t\left(\frac{1}{2 \varphi_s^2} + \frac{\gamma_1}{\varphi_s}\right)\right)_s = 0,
\end{equation}
\begin{equation} \label{ws_cl_05-2}
\Lambda^{+}_4 = e^{-t},
\qquad
\left( e^{-t} ({\varphi} + {\varphi}_t)\right)_t
    + \left(e^{-t}\left(\frac{1}{2 \varphi_s^2} + \frac{\gamma_1}{\varphi_s}\right)\right)_s = 0.
\end{equation}
Similarly, for $H^\prime = -\varphi$ and the generators $X^{-}_3$ and $X^{-}_4$ one obtains
\begin{equation} \label{ws_cl_05+1}
\Lambda^{-}_3 = \sin t,
\qquad
\left( {\varphi}\cos t - {\varphi}_t \sin t \right)_t
    - \left(
        \sin t \left(\frac{1}{2 \varphi_s^2} + \frac{\gamma_1}{\varphi_s}\right)
    \right)_s = 0,
\end{equation}
\begin{equation} \label{ws_cl_05+2}
\Lambda^{-}_4 = \cos t,
\qquad
\left( {\varphi}\sin t + {\varphi}_t \cos t \right)_t
    + \left(
        \cos t \left(\frac{1}{2 \varphi_s^2} + \frac{\gamma_1}{\varphi_s}\right)
    \right)_s = 0.
\end{equation}
To the best of authors' knowledge, conservation laws (\ref{ws_cl_05-1})--(\ref{ws_cl_05+2}) have no clear physical interpretation.

\smallskip

Finally, in case $H^\prime = 1/\varphi$, the generator~$X_5$ does not satisfy the Noether theorem.

\smallskip

One can find a more detailed discussion
of the symmetries and conservation laws of the modified shallow water equations
with different bottom topographies in~\cite{bk:KaptsovMeleshko_1D_classf[2018]}.

\begin{remark}
In case of an inclined bottom $H(z) = C_1 z + C_2$, where $z=z(t^\prime, s^\prime)$ the modified shallow
water equations in Lagrangian coordinates are reduced to the equations for a horizontal bottom $H = \text{const}$ by means of
the transformation
\begin{equation} \label{InclToFlatBottomTr}
z = \varphi + \frac{C_1}{2}  t^2,
\qquad
t^\prime = t,
\qquad
s^\prime = s.
\end{equation}
The same transformation was used in~\cite{dorodnitsyn2019shallow} for the
one--dimensional shallow water equations in Lagrangian coordinates.
\end{remark}

\subsection{Conservation laws in mass Lagrangian coordinates}

It is often possible to represent equations originally given in Eulerian or Lagrangian
coordinates in a simpler form using \emph{mass} Lagrangian coordinates.
Equations represented in mass coordinates are also often more suitable for solving problems numerically~\cite{bk:SamarskyPopov_book[1992],bk:YanenkRojd[1968]}.
Moreover, the transition from Lagrangian coordinates to mass Lagrangian coordinates is especially simple.
The mass Lagrangian coordinates are introduced with the differential form
\begin{equation} \label{dsFrom}
ds = \rho dx - u dt.
\end{equation}
Taking into account $x = \varphi(t, s)$, from the latter one derives
\begin{equation} \label{ContactTrLagr}
\varphi_t = u,
\qquad
\varphi_s = \frac{1}{\rho}.
\end{equation}
Equations~$\varphi_{ts}=\varphi_{st}$ and~(\ref{modSW_Lagr}) are brought to
the one-dimensional modified shallow water equations
in \emph{mass}~Lagrangian coordinates, namely
\begin{equation} \label{eqnsMassCoords}
\def\arraystretch{1.75}
\begin{array}{c}
\displaystyle
\left(\frac{1}{\rho}\right)_t - u_s = 0,
\\
u_t + \rho\rho_s + \gamma_1 \rho_s - H^\prime = 0,
\end{array}
\end{equation}
where $H = H(x)$.

\begin{remark}
In discretization process we will use another representation of the system introducing a new variable~\cite{dorodnitsyn2019shallow,bk:DorKapJMPSW2021,bk:DorodKapMelGN2020}
\begin{equation} \label{stateEqOrig}
p = \rho^2.
\end{equation}
Then, one can rewrite the second equation of~(\ref{eqnsMassCoords}) as
\begin{equation}
u_t + \frac{p_s}{2} + \gamma_1 \rho_s - H^\prime  = 0,
\qquad
\text{or}
\qquad
u_t + \frac{1}{2} \left(1 + \frac{\gamma_1}{\rho}\right) p_s - H^\prime
= 0.
\end{equation}
Below we will use such representation to reduce a mesh stencil for finite-difference schemes.
\end{remark}

\medskip

Consider the conservation laws possessed by system~(\ref{eqnsMassCoords}).
The first equation of~(\ref{eqnsMassCoords}) is the conservation law of mass.
By means of~(\ref{ContactTrLagr}) the conservation laws of energy~(\ref{ws_cl_01}) and momentum~(\ref{ws_cl_02}) in mass Lagrangian coordinates
are brought to
\begin{equation} \label{clm_energy}
\left(
    \frac{u^2}{2} + \frac{\rho}{2} + \gamma_1 \ln\rho - H
\right)_t
+ \left(
    u \left( \frac{\rho^2}{2} + \gamma_1 \rho \right)
\right)_s = 0,
\end{equation}
\begin{equation}
\left(
    \frac{u}{\rho}
\right)_t
+ \left(
    \rho - \frac{u^2}{2}
    + \gamma_1 \ln \rho
    - H
\right)_s = 0.
\end{equation}
The conservation laws~(\ref{ws_cl_03}) and~(\ref{ws_cl_04}) become
\begin{equation} \label{clm_momentum}
u_t
+ \left(
    \frac{\rho^2}{2} + \gamma_1 \rho
\right)_s = 0,
\end{equation}
\begin{equation} \label{clm_cmass}
\left(
    t u - x
\right)_t
+ \left(
    \frac{t\rho^2}{2} + \gamma_1 t \rho
\right)_s = 0.
\end{equation}
The conservation laws (\ref{ws_cl_05-1}) and (\ref{ws_cl_05-2}) become
\begin{equation} \label{clm_exp1}
\left(
    e^t (x - u)
\right)_t
- \left(
    e^t \left( \frac{\rho^2}{2} + \gamma_1 \rho \right)
\right)_s = 0,
\end{equation}
\begin{equation} \label{clm_exp2}
\left(
    e^{-t} (x + u)
\right)_t
+ \left(
    e^{-t} \left( \frac{\rho^2}{2} + \gamma_1 \rho \right)
\right)_s = 0,
\end{equation}
and the conservation laws (\ref{ws_cl_05+1}) and (\ref{ws_cl_05+2}) are
\begin{equation} \label{clm_sin}
\left(
    x \cos t - u \sin t
\right)_t
- \left(
    \sin t \left( \frac{\rho^2}{2} + \gamma_1 \rho \right)
\right)_s = 0,
\end{equation}
\begin{equation} \label{clm_cos}
\left(
    x \sin t + u \cos t
\right)_t
+ \left(
    \cos t \left( \frac{\rho^2}{2} + \gamma_1 \rho \right)
\right)_s = 0.
\end{equation}

\subsection{Conservation laws in Eulerian coordinates}

Denote by $D^E_t$ and $D^L_t$ Eulerian and Lagrangian total differentiation operators with respect to~$t$.
Total differentiation operators with respect to~$s$ and~$x$ are denoted by~$D_s$ and~$D_x$.
By means of~(\ref{dsFrom}), the total differentiations in Eulerian and Lagrangian coordinates are related as follows
\begin{equation}
D^L_t = D_t^E + u D_x,
\qquad
D_s = \frac{1}{\rho} D_x.
\end{equation}
The conserved quantities $(T^t, T^s)$ of conservation laws
in Lagrangian coordinates are related to the corresponding quantities $({}^eT^t, {}^eT^x)$ in Eulerian coordinates as
\begin{equation} \label{LagrEulerCLsRel}
{}^eT^t = \rho T^t,
\qquad
{}^eT^x = \rho u T^t + T^s.
\end{equation}
Formulas~(\ref{LagrEulerCLsRel}) allow one to obtain the following
Eulerian counterparts of the conservation laws~(\ref{clm_energy})--(\ref{clm_cos})
\begin{equation}
D_t^E\left(
    \rho \left(
        \frac{u^2 + \rho}{2}
        + \gamma_1 \ln \rho
        - H(x)
    \right)
\right)
+ D_x\left(
    \rho u \left(
        \frac{u^2}{2}
        + \rho
        + \gamma_1 (1 + \ln \rho)
        - H(x)
    \right)
\right) = 0,
\end{equation}
\begin{equation}
D_t^E\left(
   u
\right)
+ D_x\left(
    \frac{u^2}{2}
    + \rho
    + \gamma_1 \ln \rho
    - H(x)
\right) = 0.
\end{equation}
In case $H = \text{const}$,
\begin{equation}
D_t^E\left(
    \rho u
\right)
+ D_x\left(
    \rho u^2 + \frac{\rho^2}{2} + \gamma_1 \rho
\right) = 0.
\end{equation}
\begin{equation}
D_t^E\left(
    \rho (t u- x)
\right)
+ D_x\left(
    \rho u (t u- x) + \frac{t \rho^2}{2} + t \gamma_1 \rho
\right) = 0.
\end{equation}
In case $H = \frac{x^2}{2}$,
\begin{equation}
D_t^E\left(
    e^t \rho (x - u)
\right)
- D_x\left(
    e^t \rho \left( \frac{\rho}{2} + u^2 - x u + \gamma_1 \right)
\right) = 0.
\end{equation}
\begin{equation}
D_t^E\left(
    e^{-t} \rho (x + u)
\right)
+ D_x\left(
    e^{-t} \rho \left( \frac{\rho}{2} + u^2 + x u + \gamma_1 \right)
\right) = 0.
\end{equation}
In case $H = -\frac{x^2}{2}$,
\begin{equation}
D_t^E\left(
    \rho\left(
        x \cos t - u \sin t
    \right)
\right)
+ D_x\left(
    \rho \left(
        x u \cos t
        - \left( u^2 + \frac{\rho}{2} + \gamma_1 \right) \sin t
    \right)
\right) = 0.
\end{equation}
\begin{equation}
D_t^E\left(
    \rho\left(
        x \sin t + u \cos t
    \right)
\right)
+ D_x\left(
    \rho \left(
        x u \sin t
        + \left( u^2 + \frac{\rho}{2} + \gamma_1 \right) \cos t
    \right)
\right) = 0.
\end{equation}

\section{Discretization of the one-dimensional modified shallow water equations in Lagrangian coordinates}
\label{sec:discr}

In the present section conservative invariant finite-difference schemes
are constructed for the one-dimensional modified shallow water equations in Lagrangian
and mass Lagrangian coordinates.
The schemes early
constructed~\cite{dorodnitsyn2019shallow,bk:DorKapSW2020}  for the
shallow water equations are taken as starting point.

\medskip

Schemes in Lagrangian coordinates are further considered on a 9-point stencil
\begin{equation} \label{stencil9}
\displaystyle
(x^{n}_{m}, x^{n}_{m-1}, x^{n}_{m+1}, x^{n+1}_{m}, x^{n+1}_{m-1}, x^{n+1}_{m+1}, x^{n}_{m-1}, x^{n-1}_{m-1}, x^{n-1}_{m+1})
\equiv
(x, x_-, x_+, \hat{x}, \hat{x}_-, \hat{x}_+, \check{x}, \check{x}_-, \check{x}_+),
\end{equation}
where, in accordance with
 the notation
of~\cite{dorodnitsyn2019shallow}, we denote~$x=x^n_m=\varphi$. The
indices~$n$ and~$m$ are changed along time and space axes~$t$
and~$s$ correspondingly. The time and space steps are defined as
\begin{equation}
\def\arraystretch{1.25}
\begin{array}{c}
\displaystyle
\tau_n = \hat{\tau} = t_{n+1} - t_n = \hat{t} - t,
\qquad
\tau_{n-1} = \check{\tau} = t_n - t_{n-1} = t - \check{t},
\\
\displaystyle
h_m = h_+ = s_{m+1} - s_m = s_+ - s,
\qquad
h_{m-1} = h_- = s_{m} - s_{m-1} = s - s_-.
\end{array}
\end{equation}

\medskip

\begin{figure}[H]
  \centering
  \includegraphics[width=0.3\linewidth]{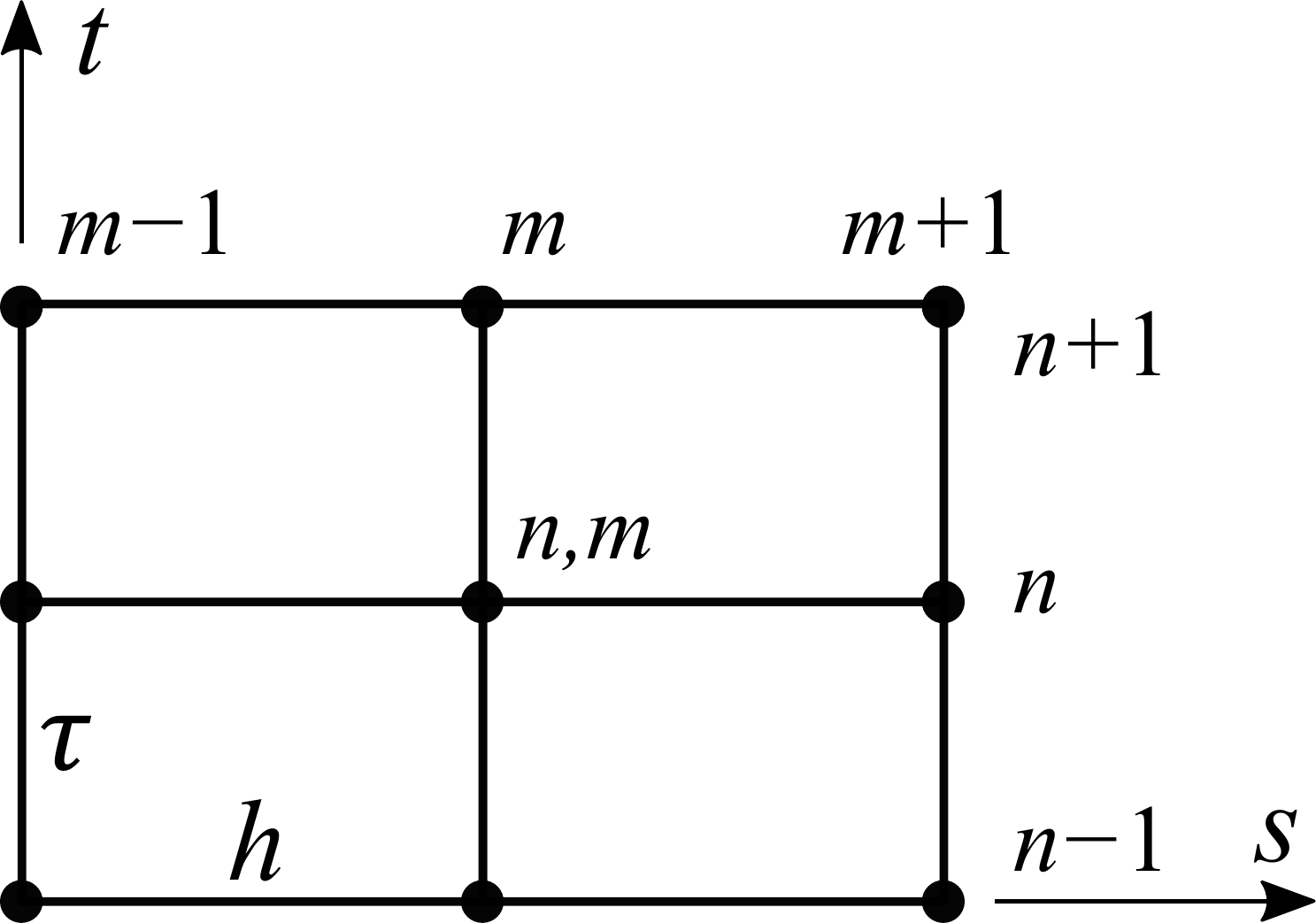}
  \caption{9-point stencil, three time layers.}
  \label{pic:stencil9}
\end{figure}

We seek for finite-difference approximations for equation~(\ref{modSW_Lagr})
defined on uniform orthogonal meshes as depicted in Figure~\ref{pic:stencil9}.

Consider the following scheme on a uniform orthogonal mesh% on stencil~(\ref{stencil9})
\begin{subequations}
\begin{equation}\label{ScmGen}
\Phi(x, x_-, x_+, \hat{x}, \hat{x}_-, \hat{x}_+, \check{x}, \check{x}_-, \check{x}_+) = 0,
\end{equation}
\begin{equation}\label{ScmGenMesh}
\displaystyle h_+ = h_- = h, \qquad
\hat{\tau} = \check{\tau} = \tau, \qquad  (\vec{\tau}, \vec{h})=0.
\end{equation}
\end{subequations}
The group generator
\begin{equation} \label{Xgen}
X = \xi^t \frac{\partial}{\partial t}
+ \xi^s \frac{\partial}{\partial s}
+ \eta \frac{\partial}{\partial x}
\end{equation}
is prolonged to the finite-difference space as follows~\cite{Dor_1,bk:Dorodnitsyn[2011]}
\begin{equation}
\displaystyle
\tilde{X} = \sum_{k,l=-\infty}^{\infty} \dtm{S^k}\dsm{S^l}(X),
\end{equation}
where the finite-difference shifts along the time and space axes are
\[
    \def\arraystretch{1.5}
    \begin{array}{c}
    \displaystyle
    \underset{\pm\tau}{S}(f(t_n, s_m, x^{n}_{m}))
        = f(t_{n \pm 1}, s_m, x^{n \pm 1}_{m}),
    \\
    \displaystyle
    \underset{\pm{s}}{S}(f(t_n, s_m, x^{n}_{m}))
        = f(t_n, s_{m\pm 1}, x^{n}_{m \pm 1}).
    \end{array}
\]

The criterion of invariance of system~(\ref{ScmGen}),~(\ref{ScmGenMesh}) is formulated as follows~\cite{bk:Dorodnitsyn[2011]}
\begin{equation} \label{SchInvCond}
\def\arraystretch{1.5}
\begin{array}{l}
\displaystyle
\tilde{X}\Phi|_{(\ref{ScmGen}), (\ref{ScmGenMesh})} = 0,
\\
\displaystyle
\tilde{X}(\hat{\tau} - \check{\tau})|_{(\ref{ScmGen}), (\ref{ScmGenMesh})} = 0,
\quad
\tilde{X}(h_+ - h_-)|_{(\ref{ScmGen}), (\ref{ScmGenMesh})} = 0,
\quad
\tilde{X}(\vec{\tau}, \vec{h})|_{(\ref{ScmGen}), (\ref{ScmGenMesh})} = 0.
\end{array}
\end{equation}
Scheme~(\ref{ScmGen}),~(\ref{ScmGenMesh}) is called \emph{invariant} if conditions~(\ref{SchInvCond}) hold.

Notice that symmetries (\ref{LagrXkern})--(\ref{SymExt3}) transform variables $t$ and $s$ independently of solution.
Therefore, the criterion of mesh invariance can be considered separately.

To preserve uniformness and orthogonality of the mesh it is needed~\cite{Dor_1,bk:Dorodnitsyn[2011]}
\begin{equation} \label{mesh_conds_uni}
  \dsp{D}\dsm{D}(\xi^s) = 0,
  \qquad
  \dtp{D}\dtm{D}(\xi^t) = 0,
\end{equation}
\begin{equation} \label{mesh_conds_ortho}
  \underset{\pm{s}}{D}(\xi^t) = -\underset{\pm{\tau}}{D}(\xi^s),
\end{equation}
where $\underset{\pm\tau}{D}$ and $\underset{\pm{s}}{D}$ are
finite-difference differentiation operators
\[
    \underset{+\tau}{D} = \frac{\underset{+\tau}{S} - 1}{t_{n+1} - t_{n}},
    \quad
    \underset{-\tau}{D} = \frac{1 - \underset{-\tau}{S}}{t_n - t_{n-1}},
    \quad
    \underset{+s}{D} = \frac{\underset{+s}{S} - 1}{s_{m+1} - s_{m}},
    \quad
    \underset{-s}{D} = \frac{1 - \underset{-s}{S}}{s_m - s_{m-1}}.
\]
One can verify that symmetries (\ref{LagrXkern})--(\ref{SymExt3}) satisfy conditions (\ref{mesh_conds_uni}) and (\ref{mesh_conds_uni}). This means that we can construct symmetry-preserving schemes for
the modified shallow water equations in Lagrangian coordinates on an invariant uniform orthogonal mesh.

\smallskip

We are seeking for conservative schemes, i.e. schemes possessing finite-difference conservation laws.
A finite-difference conservation law of scheme~(\ref{ScmGen}) is a divergent expression of the form
\begin{equation} \label{ScmCLgen}
\dtm{D}(T^t)+ \dsm{D}(T^s) = 0
\end{equation}
that vanishes on solutions of~(\ref{ScmGen}),~(\ref{ScmGenMesh}).
The quantities~$T^t$ and~$T^s$ are called density and flux.
Notice that expressions of the form~(\ref{ScmCLgen}) approximate \emph{local} differential conservation laws,
so they may not hold for discontinuous solutions such as shock waves.

\subsection{Conservative schemes for the modified shallow water equations in Lagrangian coordinates}
\label{sec:ScmsConstruction}

Consider the invariant finite-difference scheme on a uniform orthogonal mesh
for the one-dimensional shallow water equations constructed by the authors in~\cite{dorodnitsyn2019shallow}.
\begin{equation} \label{schemeSW}
\def\arraystretch{1.75}
\begin{array}{c}
\displaystyle
x_{t\check{t}}
+ \dsm{D}\left(
    \frac{1}{2\check{x}_s\hat{x}_s}
\right)
- \frac{\dtm{D}(H^h + \hat{H}^h)}{x_t + \check{x}_t}
= 0,
\\
\displaystyle
h_+ = h_- = h,
\qquad
\hat{\tau} = \check{\tau} = \tau,
\end{array}
\end{equation}
where $H^h = H^h(x)$ is some approximation for the function~$H(x)$,
$h$ and~$\tau$ are constant.

The scheme possesses the conservation law of energy
\begin{equation} \label{schemeSWenergyCL}
\displaystyle
\dtm{D}\left(
    \frac{x_t^2}{2}
    + \frac{1}{4 x_s}
    + \frac{1}{4 \hat{x}_s}
    - \frac{H^h + \hat{H}^h}{2}
\right)
+ \dsm{D}\left(
\frac{x_t^+ + \check{x}^+_t}{4\check{x}_s\hat{x}_s}
\right) = 0
\end{equation}
with the corresponding conservation law multiplier
\begin{equation} \label{LambdahEnergy}
\Lambda_1^h = \frac{x_t + \check{x}_t}{2}.
\end{equation}
The conservation law of mass
\begin{equation} \label{scmLagrCLmass}
\dtm{D}(\hat{x}_s) - \dsm{D}(x_t^+) = 0
\end{equation}
holds automatically on a uniform orthogonal mesh.

In case $H^h = \text{const}$, scheme~(\ref{schemeSW}) also possesses the conservation law
of momentum
\begin{equation}
\Lambda^h_3 = 1,
\qquad
\dtm{D}\left(x_t\right)
+ \dsm{D}\left(
    \frac{1}{2\check{x}_s\hat{x}_s}
\right)
= 0
\end{equation}
and center of mass law
\begin{equation}
\Lambda^h_4 = t,
\qquad
\dtm{D}\left(t x_t - x\right)
+ \dsm{D}\left(
    \frac{t}{2\check{x}_s\hat{x}_s}
\right)
= 0.
\end{equation}
Notice that there are additional finite-difference conservation laws exist for some specific bottom topographies~\cite{bk:DorKapJMPSW2021,dorodnitsyn2019shallow}.
We will consider them in the sequel of the paper.

\bigskip

In~\cite{bk:DorodKapMelGN2020} an invariant conservative finite-difference scheme for
the one-dimensional Green--Naghdi equations in Lagrangian coordinates was constructed by extending scheme~(\ref{schemeSW}).
A nonlinear invariant higher-order finite-difference term was added to the first equation of the scheme.
That allowed the authors to construct the scheme possessing finite-difference analogues of all the differential local conservation laws of the one-dimensional Green--Naghdi equations.
Scheme~(\ref{schemeSW}) was also used as a basis to construct
invariant conservative schemes for the two-dimensional shallow water equations in~\cite{bk:DorKapMelSW2DarXiv}.
Since the modified shallow water equations extend the shallow water equations,
it seems natural to extend scheme~(\ref{schemeSW}) in such a way that it approximates the modified shallow water equations.
Thus, here we again refer to the idea of extending the scheme previously constructed for a simpler model.

The peculiarity of such an extension is related to the fact that the conservation law of energy~(\ref{ws_cl_01})
for the modified shallow water equations includes a logarithmic term.
In contrast to the differential case, in finite differences one cannot pass from logarithmic expressions to rational expressions by means of differentiation or integration. This imposes restrictions on possible approximations for a conservative difference scheme.

%%%%%%%%%%%%%%%%%%%%%%
\medskip

Suppose that there is a finite-difference scheme that can be represented in
the form of some rational function on 9-point stencil~(\ref{stencil9}).
Then, it is natural to extend scheme~(\ref{schemeSW}) for equation~(\ref{modSW_Lagr}) as follows
\begin{equation}
\def\arraystretch{1.75}
\begin{array}{c}
\displaystyle
F_0 = x_{t\check{t}}
+ \dsm{D}\left(
    \frac{1}{2\check{x}_s\hat{x}_s}
+ \frac{\gamma_1}{\alpha_1 x_s + \alpha_2 \hat{x}_s + (1 - \alpha_1 - \alpha_2) \check{x}_s}
\right)
- \frac{\dtm{D}(H^h + \hat{H}^h)}{x_t + \check{x}_t}
= 0,
\\
\displaystyle
h_+ = h_- = h,
\qquad
\hat{\tau} = \check{\tau} = \tau,
\end{array}
\end{equation}
where $\alpha_1$ and $\alpha_2$ are some constant coefficients.

Scheme~(\ref{schemeSW}) possesses the conservation law of energy~(\ref{schemeSWenergyCL}) with the conservation law multiplier~(\ref{LambdahEnergy}).
If the extended scheme possesses an energy conservation law with multiplier~(\ref{LambdahEnergy}), then
\begin{equation} \label{badApprxExample}
\Lambda_1^h F_0 =
(x_t + \check{x}_t) \left(
    x_{t\check{t}}
    + \dsm{D}\left(
        \frac{1}{2\check{x}_s\hat{x}_s}
    + \frac{\gamma_1}{\alpha_1 x_s + \alpha_2 \hat{x}_s + (1 - \alpha_1 - \alpha_2) \check{x}_s}
    \right)
    - \frac{\dtm{D}(H^h + \hat{H}^h)}{x_t + \check{x}_t}
\right)
\end{equation}
is a finite-difference divergent expression at least for some particular values of~$\alpha_1$ and~$\alpha_2$.
To find out if this is really the case, one uses the finite-difference analogue of the direct method~\cite{bk:ChevDorKap2020,dorodnitsyn2019shallow}.
The direct method requires to consider the equation
\begin{equation} \label{direct}
\displaystyle
\mathcal{E}_x (\Lambda_1^h F_0)|_{h_- = h_+, \hat{\tau} = \check{\tau}} = 0,
\end{equation}
where $\mathcal{E}_x$ is the Euler operator on a uniform orthogonal
mesh at point~$x$~\cite{bk:Dorodnitsyn[2011]}
\begin{equation}
\displaystyle
\mathcal{E}_x = \frac{\partial}{\partial x}
- \sum_{k=-\infty}^\infty \sum_{l=-\infty}^\infty \dtm{S^k} \dsm{S^l} \left[
        \dtm{D}\left(\frac{\partial}{\partial \dtp{S^k} \dsp{S^l}(x_t)}\right)
        + \dsm{D}\left(\frac{\partial}{\partial \dtp{S^k} \dsp{S^l}(x_s)}\right)
\right]
\end{equation}
taking into account the commutation relation~(\ref{scmLagrCLmass}).

\smallskip

But, as one can verify by direct computation, there are no values~$\alpha_1$, $\alpha_2$
satisfying equation~(\ref{direct}).
Thus, approximation~(\ref{badApprxExample}) does not possess a conservation law of energy under the chosen constraints.
Despite of the fact that the considered case is quite simple,
it can be verified that more complex rational approximations for equation~(\ref{modSW_Lagr})
as well as for the integrating factor~$\Lambda^1 = x_t$ still
do not bring one closer to the goal.
The direct method is only effective enough when considering polynomial or rational expressions.
In the most general case, it is problematic to obtain the form of a scheme or a conservation law multiplier with its help.
The problem requires some additional assumptions or an ansatz to be involved.

\medskip
%%%%%%%%%%%%%%%%%%%%%%

The following considerations give a better result.

One can notice that a scheme possessing a conservation law of energy has to include some logarithmic terms.
To demonstrate this, consider the following terms of~(\ref{ws_cl_01})
\begin{equation} \label{badRel}
\left(\ln \varphi_s \right)_t - \left( \frac{\varphi_t}{\varphi_s}\right)_s
= \frac{\varphi_{st}}{\varphi_s} - \frac{\varphi_{ts}}{\varphi_s}
    - \varphi_t \left(\frac{1}{\varphi_s}\right)_s.
\end{equation}
The first and the second terms in the right hand side cancel each other out.
(This is necessary in order to write~(\ref{ws_cl_01}) as the product of~(\ref{modSW_Lagr}) and the multiplier~$\varphi_t$.)
On the contrary, in the finite-difference case this is not the true as one can see from the following
\begin{equation} \label{badApprox}
\dtm{D}\left(\ln x_s\right) - \dsm{D}\left(\frac{x_t}{x_s}\right)
= \frac{1}{\tau} \ln \frac{x_s}{\check{x}_s} - \frac{1}{x_s} x_{t\bar{s}} - x_t^- \dsm{D}\left(\frac{1}{x_s}\right).
\end{equation}
The first and the second terms of the latter expression are obviously not canceled.
This suggests that the term under~$\dsm{D}$ in the left hand side of~(\ref{badApprox}) should also include a logarithmic expression. Then, one can notice that~(\ref{badRel}) can be rewritten in the following equivalent form
\begin{equation} \label{goodRel}
\left(\ln \varphi_s \right)_t - \left( \frac{\varphi_t}{\varphi_{st}} \, (\ln \varphi_s)_t \right)_s.
\end{equation}
The latter can be used as an ansatz for construction of conservative schemes.
For example, one can construct the following finite-difference analogue of~(\ref{goodRel})
\begin{equation} \label{approxResult0}
\dtm{D}\left(\ln x_s\right) - \dsm{D}\left( x_t^+ \, \frac{\dtm{D} (\ln x_s)}{x_{st}} \right)
=
\frac{1}{\tau} \ln \frac{x_s}{\check{x}_s}
- \frac{x_{ts}}{\tau x_{st}} \ln \frac{x_s}{\check{x}_s}
- x_t \dsm{D}\left( \frac{1}{\tau x_{st}} \ln \frac{x_s}{\check{x}_s} \right)
= - x_t \dsm{D}\left( \frac{1}{\tau x_{st}} \ln \frac{x_s}{\check{x}_s} \right)
\end{equation}
which is an approximation of a desired form.
The resulting expression~(\ref{approxResult0}) corresponds to the conservation law multiplier~$x_t$.
For the multiplier~(\ref{LambdahEnergy}) some minor changes required,
namely one should consider the expression
\begin{equation} \label{approxResult1}
\dtm{D}\left(
    \ln (x_s \hat{x}_s )
\right)
- \dsm{D}\left(
        \frac{x_t^+ + \check{x}^+_t}{\tau (x_s + \check{x}_s)_t} \, \ln \frac{\hat{x}_s}{\check{x}_s}
\right)
= - \Lambda^h_1 \dsm{D}\left( \frac{2}{\tau (x_s + \check{x}_s)_t} \, \ln \frac{\hat{x}_s}{\check{x}_s} \right)
= \Lambda^h_1 \dsm{D}\left( \frac{2}{\check{x}_s - \hat{x}_s} \, \ln \frac{\hat{x}_s}{\check{x}_s} \right).
\end{equation}
instead of~(\ref{approxResult0}).
The latter can be obtained by means of algebraic transformations or, more systematically, with the help of the direct method.

\medskip

Based on the above, one extends scheme~(\ref{schemeSW}) to
the following invariant conservative scheme for the one-dimensional modified shallow water equations~(\ref{modSW_Lagr})
\begin{equation} \label{scheme}
\def\arraystretch{2}
\begin{array}{c}
\displaystyle
x_{t\check{t}}
+ \dsm{D}\left(
    \frac{1}{2\check{x}_s\hat{x}_s}
\right)
+ \dsm{D}\left(
    \frac{\gamma_1}{\hat{x}_s - \check{x}_s} \, \ln \frac{\hat{x}_s}{\check{x}_s}
\right)
- \frac{\dtm{D}(H^h + \hat{H}^h)}{x_t + \check{x}_t}
= 0,
\\
\displaystyle
h_+ = h_- = h,
\qquad
\hat{\tau} = \check{\tau} = \tau.
\end{array}
\end{equation}
This scheme is defined on 9-point stencil~(\ref{stencil9}).
It approximates~(\ref{modSW_Lagr}) up to~$O(\tau^2 + h^2)$.

\medskip

Scheme~(\ref{scheme}) is an invariant one.
In case $H^h=H^h(x)$, it admits the generators $X_1$ and $X_2$.
In case $H^h=\text{const}$, it also admits the generators $X_3$, $X_4$ and $X_5$.

\medskip

Scheme~(\ref{scheme}) possesses conservation law of mass~(\ref{scmLagrCLmass})
and the conservation law of energy (with the conservation law multiplier
$\Lambda^h_1 = \frac{1}{2}(x_t + \check{x}_t)$)
\begin{equation} \label{SchemeCLenergyH}
\displaystyle
\dtm{D}\left(
    \frac{x_t^2}{2}
    + \frac{1}{4 x_s}
    + \frac{1}{4 \hat{x}_s}
    - \frac{\gamma_1}{2} \ln (x_s \hat{x}_s )
    - \frac{H^h + \hat{H}^h}{2}
\right)
+ \dsm{D}\left(
    \frac{x_t^+ + \check{x}^+_t}{2}
    \left(
        \frac{1}{2\check{x}_s\hat{x}_s}
        + \frac{\gamma_1}{\hat{x}_s - \check{x}_s} \, \ln \frac{\hat{x}_s}{\check{x}_s}
    \right)
\right) = 0.
\end{equation}

Notice that there is no finite-difference analogue of the conservation law~(\ref{ws_cl_02}).
As alternative, one can construct schemes that conserve momentum but do not possess the energy conservation law.
See discussion in~\cite{dorodnitsyn2019shallow}.
For a horizontal bottom topography one can also consider the conservation law of momentum~(\ref{scmLagrmom}).

\subsubsection{Case of a horizontal bottom}

In case $H^h = \text{const}$, there are two additional conservation laws, \emph{momentum} and \emph{center-of-mass law},
with the conservation law multipliers $1$  and $t$.
\begin{equation} \label{scmLagrmom}
\dtm{D}\left(x_t\right)
+ \dsm{D}\left(
    \frac{1}{2\check{x}_s\hat{x}_s}
    + \frac{\gamma_1}{\hat{x}_s - \check{x}_s} \, \ln \frac{\hat{x}_s}{\check{x}_s}
\right)
= 0,
\end{equation}
\begin{equation} \label{scmLagrcmass}
\dtm{D}\left(t x_t - x\right)
+ \dsm{D}\left(
    \frac{t}{2\check{x}_s\hat{x}_s}
    + \frac{\gamma_1 t}{\hat{x}_s - \check{x}_s} \, \ln \frac{\hat{x}_s}{\check{x}_s}
\right)
= 0.
\end{equation}

\subsubsection{Case of an inclined bottom}

In case~$H(z) = C_1 z + C_2$ (i.e., $H^\prime = C_1$), where $z=z(t^\prime, s^\prime)$, the corresponding scheme
\begin{equation}
\def\arraystretch{1.75}
\begin{array}{c}
\displaystyle
z_{t^\prime\check{t}^\prime}
+ \underset{-s^\prime}{D}\left(
    \frac{1}{2\check{z}_{s^\prime}\hat{z}_{s^\prime}}
\right)
+ \underset{-s^\prime}{D}\left(
    \frac{\gamma_1}{\hat{z}_{s^\prime} - \check{z}_{s^\prime}} \, \ln \frac{\hat{z}_{s^\prime}}{\check{z}_{s^\prime}}
\right) - C_1
= 0,
\\
{s^\prime}_+ - {s^\prime} = {s^\prime} - {s^\prime}_- = h^\prime,
\\
\hat{t}^\prime - t^\prime = {t^\prime} - \check{t}^\prime = \tau^\prime
\end{array}
\end{equation}
can be transformed into the form~(\ref{scheme})
by means of the following finite-difference analogue of~(\ref{InclToFlatBottomTr})
\begin{equation} \label{InclToFlatBottomTrDelta}
z = x + \frac{C_1}{2} t \hat{t},
\qquad
t' = t,
\qquad
s' = s.
\end{equation}
Thus, in Lagrangian coordinates the case of an inclined bottom is reduced to the case of a horizontal bottom~(see also~\cite{dorodnitsyn2019shallow}).
%Notice that in Eulerian coordinates there
%is an analog of transformation (1) for the shallow water equations~\cite{bk:ChirkunovPikmullina[2014]}.

\subsubsection{Case of parabolic bottoms}

In the paper~\cite{bk:DorKapJMPSW2021} the authors have constructed the schemes for the parabolic bottoms $H(x) = \pm\frac{x^2}{2}$. By analogy to the previous case, they are extended as follows.

\smallskip

In case $H(x) = \frac{x^2}{2}$, the approximation
\begin{equation}
H^h = \frac{\cosh \tau - 1}{\tau^2} x\hat{x} = \frac{x^2}{2} + O(\tau)
\end{equation}
for the function $H(x)$ was found with the help of the direct method.
It was found that the latter approximation is necessary for the scheme
to possess the conservation laws with the conservation law multipliers~$e^{\pm{t}}$.
(See~\cite{bk:DorKapJMPSW2021} for details.)
Thus, the extended scheme for the modified shallow water equations is
\begin{equation} \label{scheme_sq_lagr}
  \displaystyle
  \def\arraystretch{1.75}
  \begin{array}{c}
        \displaystyle
        x_{t\check{t}}
      + \dsm{D}\left(
            \frac{1}{2\check{x}_s\hat{x}_s}
            + \frac{\gamma_1}{\hat{x}_s - \check{x}_s} \, \ln \frac{\hat{x}_s}{\check{x}_s}
      \right) - \frac{2 (\cosh \tau - 1)}{\tau^2} x
      = 0,
      \\
      \displaystyle
      \qquad
      h_+ = h_- = h,
    \qquad
    \hat{\tau} = \check{\tau} = \tau.
  \end{array}
\end{equation}
The conservation law of energy for the latter scheme is obtained by formula~(\ref{SchemeCLenergyH})
\begin{multline} \label{SchemeCLenergyHParabolic}
\Lambda_1^h = \frac{x_t + \check{x}_t}{2},
\quad
\dtm{D}\left(
    \frac{x_t^2}{2}
    + \frac{1}{4 x_s}
    + \frac{1}{4 \hat{x}_s}
    - \frac{\gamma_1}{2} \ln (x_s \hat{x}_s )
    -  \frac{2 (\cosh \tau - 1)}{\tau^2}  x\hat{x}
\right) \\
+ \dsm{D}\left(
    \frac{x_t^+ + \check{x}^+_t}{2}
    \left(
        \frac{1}{2\check{x}_s\hat{x}_s}
        +\frac{\gamma_1}{\hat{x}_s - \check{x}_s} \, \ln \frac{\hat{x}_s}{\check{x}_s}
    \right)
\right) = 0.
\end{multline}
The additional finite-difference conservation laws and the multipliers
which correspond to~(\ref{ws_cl_05-1}) and~(\ref{ws_cl_05-2}) are the following
\begin{equation} \label{scmLagrExp1}
(\Lambda_3^{+})^h = e^t,
\quad
{\dtm{D}}\left(
    x \frac{e^{\hat{t}} - e^t}{\tau}
    - e^t {x}_t
\right)
    -\dsm{D}\left(
        e^t
        \left(
            \frac{1}{2\check{x}_s\hat{x}_s}
            + \frac{\gamma_1}{\hat{x}_s - \check{x}_s}\, \ln \frac{\hat{x}_s}{\check{x}_s}
        \right)
    \right)
    = 0,
\end{equation}
\begin{equation} \label{scmLagrExp2}
(\Lambda_4^{+})^h = e^{-t},
\quad
{\dtm{D}}\left(
    x \frac{e^{-t} - e^{-\hat{t}}}{\tau}
    - e^{-t} {x}_t
\right)
    +\dsm{D}\left(
        e^{-t}
        \left(
            \frac{1}{2\check{x}_s\hat{x}_s}
            + \frac{\gamma_1}{\hat{x}_s - \check{x}_s} \, \ln \frac{\hat{x}_s}{\check{x}_s}
        \right)
    \right)
    = 0.
\end{equation}

\medskip

In case $H = -\frac{x^2}{2}$, the scheme is extended the similar way, and one derives
\begin{equation} \label{scheme_sq_lagr2}
  \displaystyle
  \def\arraystretch{1.75}
  \begin{array}{c}
        \displaystyle
        x_{t\check{t}}
      + \dsm{D}\left(
            \frac{1}{2\check{x}_s\hat{x}_s}
            + \frac{\gamma_1}{\hat{x}_s - \check{x}_s} \, \ln \frac{\hat{x}_s}{\check{x}_s}
      \right)
       - \frac{2 (\cos \tau - 1)}{\tau^2} x
      = 0,
      \\
      \displaystyle
      h_+ = h_- = h,
    \qquad
    \hat{\tau} = \check{\tau} = \tau.
  \end{array}
\end{equation}
The conservation law of energy and two additional conservation laws which correspond to~(\ref{ws_cl_05+1}) and~(\ref{ws_cl_05+2}) are
\begin{multline}
\Lambda_1^h = \frac{x_t + \check{x}_t}{2},
\quad
\dtm{D}\left(
    \frac{x_t^2}{2}
    + \frac{1}{4 x_s}
    + \frac{1}{4 \hat{x}_s}
    - \gamma_1 \ln (x_s \hat{x}_s )
    - \frac{2 (\cos \tau - 1)}{\tau^2} x\hat{x}
\right) \\
+ \dsm{D}\left(
    \frac{x_t^+ + \check{x}^+_t}{2}
    \left(
        \frac{1}{2\check{x}_s\hat{x}_s}
        +\frac{\gamma_1}{\hat{x}_s - \check{x}_s} \, \ln \frac{\hat{x}_s}{\check{x}_s}
    \right)
\right) = 0,
\end{multline}
\begin{equation} \label{scmLagrCos}
(\Lambda_3^{-})^h = \cos t,
\quad
{\dtm{D}}\left(
    {x}_t \cos t
    - x \frac{\cos \hat{t} - \cos t}{\tau}
\right)
    +\dsm{D}\left(
        \cos t
        \left(
            \frac{1}{2\check{x}_s\hat{x}_s}
            + \frac{\gamma_1}{\hat{x}_s - \check{x}_s}  \, \ln \frac{\hat{x}_s}{\check{x}_s}
        \right)
    \right) = 0.
\end{equation}
\begin{equation} \label{scmLagrSin}
(\Lambda_4^{-})^h = \sin t,
\quad
{\dtm{D}}\left(
    {x}_t \sin t
    - x \frac{\sin \hat{t} - \sin t}{\tau}
\right)
    +\dsm{D}\left(
        \sin t
        \left(
            \frac{1}{2\check{x}_s\hat{x}_s}
            + \frac{\gamma_1}{\hat{x}_s - \check{x}_s} \, \ln \frac{\hat{x}_s}{\check{x}_s}
        \right)
    \right) = 0.
\end{equation}

\subsubsection{Comparison with a naive approximation}

In Section~\ref{sec:numimpl} devoted to the numerical implementation of schemes,
the analysis of energy preservation by various schemes will be performed.
The present section provides some preliminary ideas.

\medskip

Consider the following naive scheme in case of the horizontal bottom $H=0$
\begin{equation} \label{naive}
\def\arraystretch{1.75}
\begin{array}{c}
\displaystyle
x_{t\check{t}}
+ \dsm{D}\left(
    \frac{1}{2\check{x}_s\hat{x}_s}
    + \frac{\gamma_1}{x_s}
\right)
= 0,
\\
\displaystyle
h_+ = h_- = h,
\qquad
\hat{\tau} = \check{\tau} = \tau.
\end{array}
\end{equation}
This scheme is invariant and approximates the modified shallow water equations up to~$O(\tau^2 + h^2)$.

Scheme~(\ref{naive}) is constructed in obvious way by extending scheme~(\ref{schemeSW}) with the rational term
\[
\dsm{D}(\gamma_1/x_s).
\]
As it was shown in the beginning of Section~\ref{sec:ScmsConstruction},
such a scheme cannot possess a conservation law of energy corresponding to the conservation law multiplier~(\ref{LambdahEnergy}). One can verify the similar way that more general multipliers of polynomial
or even rational form do not lead to a conservation law of energy.
Thus, either scheme~(\ref{naive}) does not possess a conservation law of energy,
or it possesses a conservation law of energy of a very complicated form.
% and there is no way to construct such a conservation law using methods known to the authors.

\medskip

In order to compare how the energy is preserved by schemes~(\ref{scheme}) and~(\ref{naive}),
we will use two approaches.

\medskip

\begin{enumerate}

\item
Assuming that scheme~(\ref{naive}) does not possess an energy conservation law in principle,
we construct a quite reasonable approximation of it, based on the known conservation laws~(\ref{schemeSWenergyCL})
and~(\ref{SchemeCLenergyH}) for schemes~(\ref{schemeSW}) and~(\ref{scheme}).
Based on~(\ref{schemeSWenergyCL}), we consider the following approximation for the conservation law of energy
\begin{multline} \label{deltaEpsExpr0}
\frac{x_t + \check{x}_t}{2}\left(
    x_{t\check{t}}
    + \dsm{D}\left(
        \frac{1}{2\check{x}_s\hat{x}_s}
        + \frac{\gamma_1}{x_s}
    \right)
\right)
=
\dtm{D}\left(
    \frac{x_t^2}{2}
    + \frac{1}{4 x_s}
    + \frac{1}{4 \hat{x}_s}
\right)
\\
+ \dsm{D}\left(
\frac{x_t^+ + \check{x}^+_t}{4\check{x}_s\hat{x}_s}
\right)
+ \gamma_1 \frac{x_t + \check{x}_t}{2} \dsm{D}\left(\frac{1}{x_s}\right).
\end{multline}
The term at $\gamma_1$ is not a divergent expression.
%\todo{In order (91) to approximate diff CL we add log. Then, the first part approx CL of enertgy. We add logarithmic terms.. We have nondivergent rest}
%Taking~(\ref{SchemeCLenergyH}) into account, by algebraic transformations we modify the latter
%approximation as follows
In order for expression~(\ref{deltaEpsExpr0}) to approximate conservation law~(\ref{ws_cl_01}),
a logarithmic term should be added to the density of the conservation law.
It cannot arise from rational expressions, and therefore we introduce it artificially.
There are various ways to do this, and we introduce it in the same form as in conservation law~(\ref{SchemeCLenergyH})
for scheme~(\ref{scheme}). To do this, we add and subtract new divergent terms to~(\ref{deltaEpsExpr0}).
Then, rearranging the terms, we derive
\begin{multline} \label{deltaEpsExpr}
\dtm{D}\left(
    \frac{x_t^2}{2}
    + \frac{1}{4 x_s}
    + \frac{1}{4 \hat{x}_s}
    - \frac{\gamma_1}{2} \ln x_s \hat{x}_s
\right)
+ \dsm{D}\left(
\frac{x_t^+ + \check{x}^+_t}{4\check{x}_s\hat{x}_s}
+ \gamma_1\frac{x_t^+ + \check{x}^+_t}{2 x_s}
\right)
\\
- \gamma_1 \left\{
    \frac{x_t + \check{x}_t}{2} \frac{x_{s\bar{s}}}{x_s x_{\bar{s}}}
    + \dsm{D}\left(\frac{x_t^+ + \check{x}^+_t}{2 x_s}\right)
    - \frac{1}{2} \dtm{D}(\ln x_s \hat{x}_s)
\right\}.
\end{multline}
The terms in the first two brackets form a reasonable approximation for the conservation law~(\ref{ws_cl_01}).
Indeed, with infinite mesh refinement, they vanish on the solutions of the scheme.
The expression in the curly braces tends to zero in continuous limit,
so that one can consider it as the energy preservation error~$\delta\varepsilon$.
The Taylor series expansion gives the estimation
\begin{equation} \label{deltaEpsSeries}
\displaystyle
%\delta\varepsilon = \tau \frac{\gamma_1}{2}\left(\frac{\varphi_{ts}}{\varphi_s} \right)_t + O(h\tau + \tau^2).
%\delta\varepsilon = \frac{3 \varphi_s \varphi_{ts} \varphi_{tts} - 2 \varphi_{ts}^3}{6 \varphi_s^3} \, \gamma_1 \tau^2  + O(h\tau^2).
\delta\varepsilon = \left( \left(\frac{1}{2 \varphi_s}\right)_t \varphi_{tts}
    + \frac{((\ln \varphi_s)_t)^3}{3} \right)  \gamma_1 \tau^2  + O(h\tau^2).
\end{equation}
Notice that representation~(\ref{deltaEpsExpr}) of~(\ref{deltaEpsExpr0}) is not unique.
We have chosen an appropriate approximation based on~(\ref{schemeSWenergyCL}) analyzing schemes~(\ref{naive}) and~(\ref{scheme}).
It can be shown that other ap\-prox\-i\-ma\-tions give similar results.

\item
Another approach is to measure the total energy evolution in time without regard to a particular scheme. %~\cite{bk:DorKapJMPSW2021}. %,bk:Bihlo_numeric[2012]}.
In Lagrangian coordinates, we consider the following sum
\begin{equation}
\mathcal{H}(n) = \frac{h}{2} \sum_{(i)} \left[\left(\frac{x^{n+1}_i-x^n_i}{\tau}\right)^2
    + \frac{h}{x^n_{i+1}-x^n_i} - 2\gamma_1 \ln \frac{(x^n_{i+1}-x^n_i)}{h} \right]
\end{equation}
whose value gives the total energy in the computational domain.
Its value should tend to constant in the continuous limit~\cite{bk:Dorod_Hamilt[2011]}.
$\mathcal{H}(n)$ corresponds to the total energy for the modified shallow water equations
\begin{equation}
 \tilde{\mathcal{H}}(t,s) = \int \left(\frac{\varphi_t^2}{2}
    + \frac{1}{2 \varphi_s}
    - \gamma_1 \ln \varphi_s\right) d\varphi
\end{equation}
at time $t=n\tau$.
Recall that the function under integral is the sum of kinetic and potential
energy.
To estimate the change in total energy over time, we will consider the relative error
\begin{equation} \label{ErrRelativeFormula}
e_R(n) = \frac{|\mathcal{H}(n) - \mathcal{H}(0)|}{|\mathcal{H}(0)|}
\end{equation}
for the schemes under comparison.
\end{enumerate}

The advantage of the second approach is that we do not make any assumptions about the existence of conservation laws for the schemes. On the other hand, the first approach allows one to obtain a reasonable estimate of energy preservation
for a fixed moment of time under certain assumptions about the form of the finite-difference conservation laws.

\subsection{Conservative schemes in mass Lagrangian coordinates}

%Scheme~(\ref{scheme}) is defined on three time layers.
%The statement of initial and boundary value problems for three-layer schemes
%can be quite complicated compared to two-layer schemes.
%It turns out that the number of time layers can be decreased by switching from Lagrangian coordinates to mass Lagrangian coordinates~\cite{dorodnitsyn2019shallow,bk:DorodKapMelGN2020}.
%However, a straightforward transition to mass coordinates using simple approximations of~(\ref{ContactTrLagr})
A straightforward transition to mass Lagrangian coordinates using simple approximations of~(\ref{ContactTrLagr})
leads to three-layer schemes. Now we demonstrate that by choosing an \emph{appropriate} approximation for~(\ref{ContactTrLagr}) and~(\ref{stateEqOrig})
one can rewrite the scheme in mass Lagrangian coordinates on \emph{two} time layers.
Following~\cite{dorodnitsyn2019shallow}, we choose the approximation
\begin{equation} \label{ContactTrLagrApprx}
  \check{x}_s + x_s = \frac{2}{\check{\rho}},
  \qquad
  x_t = u
\end{equation}
for equation~(\ref{ContactTrLagr}), and the implicit approximation
\begin{equation} \label{stateEqOrigApprx}
    \displaystyle
    \frac{1}{\sqrt{\check{p}}} + \frac{1}{\sqrt{p}} = \frac{2}{\check{\rho}}
\end{equation}
for~(\ref{stateEqOrig}).
Notice that by virtue of~(\ref{ContactTrLagrApprx}) the latter can be considered equivalent to~$p = 1/x_s^2$.
It is also important to note that the change of variables does not affect the independent variables~$t$ and~$s$,
so the uniform orthogonal mesh remains invariant.

\begin{figure}[H]
  \centering
  \includegraphics[width=0.4\linewidth]{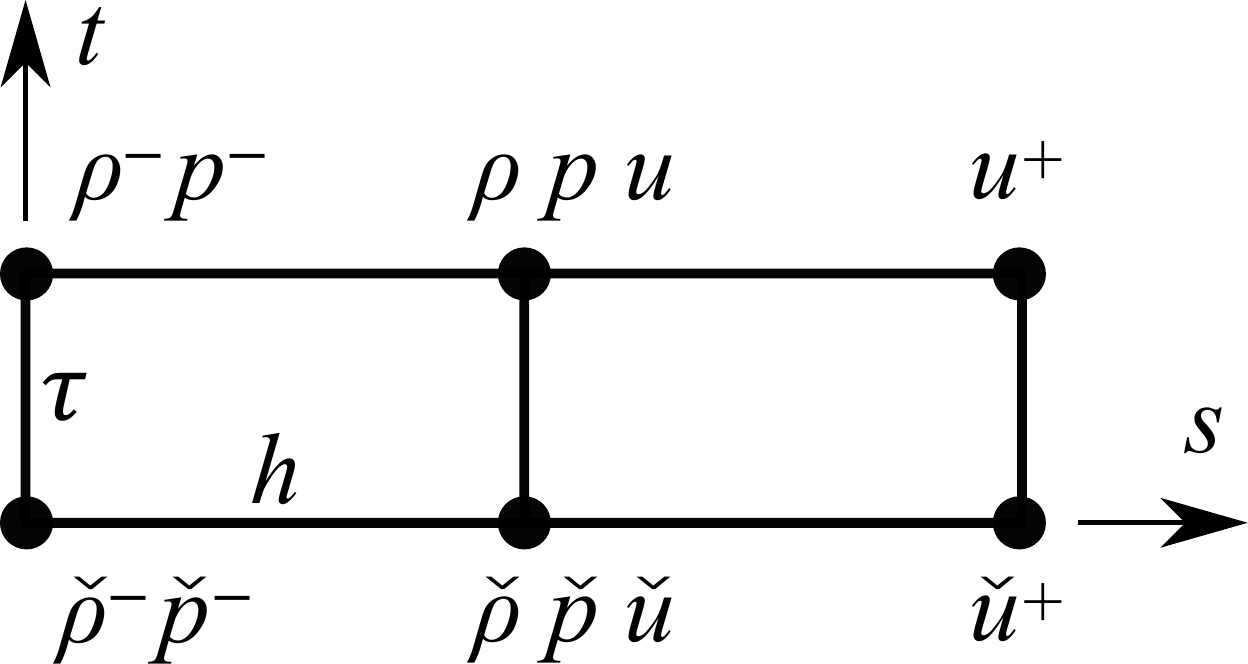}
  \caption{6-point stencil, two time layers.}
  \label{pic:stencil6}
\end{figure}

Then, scheme~(\ref{scheme})
can be expressed in mass Lagrangian coordinates on two time layers as\footnote{To derive the first equation of the scheme one gets the sum of~(\ref{scmLagrCLmass}) and the shifted one, i.e.,
\[\dtm{D}(\hat{x}_s + x_s) - \dsm{D}(x_t^+ + \check{x}_t^+) = 0.\]}
\begin{equation} \label{schemeLagrMass}
    \def\arraystretch{1.85}
    \begin{array}{c}
    \displaystyle
    \dtm{D}\left( \frac{1}{\rho} \right)
    - \dsm{D}\left(
        \frac{u^+ + \check{u}^+}{2}
    \right) = 0,
    \\
    \displaystyle
      \dtm{D}(u) + \dsm{D}\left(
        Q
    \right)
    - \frac{\dtm{D}(H^h + \hat{H}^h)}{u + \check{u}}
    = 0,
    \\
    \displaystyle
    x_t = u,
    \qquad
    \check{x}_s + x_s = \frac{2}{\check{\rho}},
    \qquad
    \frac{1}{\sqrt{\check{p}}} + \frac{1}{\sqrt{p}} = \frac{2}{\check{\rho}},
    \\
    \displaystyle
    \hat{\tau} = \check{\tau} = \tau,
    \qquad
    h_+ = h_- = \ h,
    \qquad
    (\vec{\tau},\vec{h}) = 0,
    \end{array}
\end{equation}
where
\begin{equation} \label{flux_Q}
\displaystyle
  Q = \frac{1}{2}\left[
        \frac{4}{\rho \check{\rho}}
        - \frac{2}{\sqrt{p}}\left( \frac{1}{\rho} + \frac{1}{\check{\rho}} \right)
        + \frac{1}{p}
    \right]^{-1}
     - \frac{\gamma_1 \rho\check{\rho}}{\rho - \check{\rho}} \, \ln \left[
        \sqrt{\check{p}} \left(
            \frac{2}{\rho} - \frac{1}{\sqrt{p}}
        \right)
     \right]
  = \frac{\rho^2}{2} + \gamma_1 \rho + O(\tau).
\end{equation}
Thus, scheme~(\ref{schemeLagrMass}) is defined on \emph{two} time layers (see Figure~\ref{pic:stencil6}) by including additional equations~(\ref{ContactTrLagrApprx}) and~(\ref{stateEqOrigApprx}) into the system.
Notice that for the first time such an approach was proposed in~\cite{bk:Korobitsyn_scheme[1989]}
for constructing schemes possessing an extended set of conservation laws.

\bigskip

By means of~(\ref{ContactTrLagrApprx}) and (\ref{stateEqOrigApprx})
one obtains the following finite-difference conservation laws for scheme~(\ref{schemeLagrMass}).

The conservation law of mass is the first equation of scheme~(\ref{schemeLagrMass}), i.e.,
\begin{equation}
\displaystyle
    \dtm{D}\left( \frac{1}{\rho} \right)
    - \dsm{D}\left(
        \frac{u^+ + \check{u}^+}{2}
    \right) = 0.
\end{equation}

For an arbitrary $H(x)$, the conservation law of energy~(\ref{clm_energy}) is brought to
\begin{equation} \label{scmLagrMassEnergyCL}
\dtm{D}\left(
    \frac{u^2}{2} + \frac{1}{2}\,\frac{p}{\rho - 2\sqrt{p}}
    - \gamma_1 \ln \left( \frac{2}{\rho\sqrt{p}} - \frac{1}{p} \right)
    - \frac{H^h + \hat{H}^h}{2}
\right)
+ \dsm{D}\left(
    \frac{u^+ + \check{u}^+}{2}\, Q
\right) = 0.
\end{equation}
The conservation laws of energy for specific bottom topographies are
obtained from~(\ref{scmLagrMassEnergyCL}) in quiet straightforward way so they are not presented here.

The conservation law of momentum (\ref{scmLagrmom}) and the center of mass law~(\ref{scmLagrcmass})
for the case~$H(x) = \text{const}$ in mass Lagrangian coordinates are
% іїTЂіїі-ііTЃT‚ії і-і-і-іііїT‚іііїії
\begin{equation}
\dtm{D}\left(u\right) + \dsm{D}\left(Q\right) = 0,
\end{equation}
\begin{equation}
\dtm{D}\left(t u - x\right) + \dsm{D}\left(t Q\right) = 0,
\end{equation}

The additional conservation laws (\ref{scmLagrExp1}) and (\ref{scmLagrExp2})
for a parabolic bottom profile $H(x)=\frac{x^2}{2}$ become
\begin{equation}
{\dtm{D}}\left(
    x \frac{e^{\hat{t}} - e^t}{\tau}
    - e^t u
\right)
    -\dsm{D}\left(
        Q e^t
    \right)
    = 0,
\end{equation}
\begin{equation}
{\dtm{D}}\left(
    x \frac{e^{-t} - e^{-\hat{t}}}{\tau}
    - e^{-t} u
\right)
    +\dsm{D}\left(
        Q e^{-t}
    \right)
    = 0.
\end{equation}
Similarly, the additional conservation laws (\ref{scmLagrCos}) and (\ref{scmLagrSin})
for $H(x)=-\frac{x^2}{2}$
in mass coordinates are
\begin{equation}
{\dtm{D}}\left(
    u \cos t
    - x \frac{\cos \hat{t} - \cos t}{\tau}
\right)
    +\dsm{D}\left(
        Q \cos t
    \right) = 0.
\end{equation}
\begin{equation}
{\dtm{D}}\left(
    u \sin t
    - x \frac{\sin \hat{t} - \sin t}{\tau}
\right)
    +\dsm{D}\left(
        Q \sin t
    \right) = 0.
\end{equation}

\section{Numerical implementation of the constructed schemes}
\label{sec:numimpl}

In the present section the numerical implementation of conservative
scheme~(\ref{scheme}) and ``naive'' scheme~(\ref{naive}) are
considered on two test problems. The first problem is the dam-break
problem over a parabolic bottom. For the one-dimensional shallow water
equations in Lagrangian coordinates it was considered in~\cite{bk:DorKapJMPSW2021}
where it was numerically implemented for scheme~(\ref{schemeSW}).
The second problem is the collapse of a fluid column above an inclined
bottom which was considered in~\cite{bk:KarPetSla09} for the modified
shallow water equations in Eulerian coordinates.

As explained in Section~\ref{sec:diffeqns}, for a given initial height $\rho_{0}(\xi)$
of the fluid over the bottom~$H$, one can derive the function $\alpha(s)$
by solving the Cauchy problem
\begin{equation}
\rho_{0}(\alpha(s))=\frac{1}{\alpha^{\prime}(s)},\qquad\alpha(0)=0.\label{CauchyProblem}
\end{equation}
Solving the equation $x=\varphi(t,\alpha(s))$ with respect to $s=A(t,x)$,
and using the identity
\[
x-\varphi(t,\alpha(A(t,x)))=0,
\]
one obtains that
\[
A_{x}(t,x)=\rho(\alpha(A(t,x)),t).
\]
Hence, the initial distribution~$A(t_{0},x)$ becomes
\begin{equation}
{\displaystyle A(t_{0},x)=\int_{0}^{x}\rho_{0}(\xi)d\xi.}\label{massLagrInt}
\end{equation}
For more detailed discussion on calculations in Lagrangian coordinates
see~\cite{bk:DorodKapMelGN2020,bk:DorKapJMPSW2021}.

\bigskip

To implement schemes (\ref{scheme}) and (\ref{naive}) numerically,
we represent them in the following form
\begin{equation} \label{toLinearize}
\hat{x} - 2x + \check{x}
    + \frac{h\tau^2}{2} \frac{
        (\hat{x} - \hat{x}_-)(\check{x} - \check{x}_-)
        -
        (\hat{x}_+ - \hat{x})(\check{x}_+ - \check{x})
    }
    {
        (\hat{x} - \hat{x}_-)(\check{x} - \check{x}_-)
        (\hat{x}_+ - \hat{x})(\check{x}_+ - \check{x})
    }
    + \tau^2\dsm{D}(\gamma_1\Gamma^h)
    + \tau^2 (H^\prime)^h(x) = 0,
\end{equation}
\[
\displaystyle
h_+ = h_- = h,
\qquad
\hat{\tau} = \check{\tau} = \tau,
\]
where $(H^\prime)^h$ is a chosen approximation for the derivative~$H^\prime(x)$,
and $\Gamma^h$ is an approximation for the term $1/\varphi_s$.

Then, we linearize equation~(\ref{toLinearize}) with respect to the solution on the upper time layer~\cite{bk:DorKapJMPSW2021}
\begin{multline} \label{tridiag}
    \frac{h^2\tau^2}{\Delta_1}(x^{n-1}_m - x^{n-1}_{m-1}) \, x^{(j+1)}_{m-1}
    - \left(1 + \frac{h^2\tau^2}{\Delta_1}(x^{n-1}_{m+1} - x^{n-1}_{m-1}) \right) x^{(j+1)}_{m}
    \\
    + \frac{h^2\tau^2}{\Delta_1}(x^{n-1}_{m+1} - x^{n-1}_m) \, x^{(j+1)}_{m+1}
    = 2x^n_m - x^{n-1}_m
        -\tau^2\dsm{D}(\gamma_1\Gamma^h)
        -\tau^2(H^\prime)^h(x^n_m),
\end{multline}
where the indices ${}^{(j)}$ denote the number of iteration,
\[
\Delta_1 = 2\,(x^{(j)}_m - x^{(j)}_{m-1})
        (x^{(j)}_{m+1} - x^{(j)}_m)
        (x^{n-1}_m - x^{n-1}_{m-1})
        (x^{n-1}_{m+1} - x^{n-1}_m),
\]
\[
    n = 2, 3, \dots,
    \qquad
    m = 2, 3, \dots, \lfloor L/h \rfloor - 1.
\]

The advantage of representation~(\ref{tridiag}) is that
on each iteration it can be solved with the help of tridiagonal matrix algorithm~(one can find description and the stability conditions of this well-known algorithm, e.~g., in~\cite{bk:Samarskii2001theory}).

\smallskip

For scheme~(\ref{naive}), $\Gamma^h = 1/x_s = h/(x^n_{m+1} - x^n_m)$, and for scheme~(\ref{scheme}) it is
\begin{equation} \label{hardTerm}
\Gamma^h = \frac{1}{\hat{x}_s - \check{x}_s}  \, \ln \frac{\hat{x}_s}{\check{x}_s}
= \frac{h}{x^{n+1}_{m+1} - x^{n+1}_{m} - x^{n-1}_{m+1}
    +x^{n-1}_{m}} \, \ln \frac{x^{n+1}_{m+1} - x^{n+1}_{m}}{x^{n-1}_{m+1} - x^{n-1}_{m}}.
\end{equation}
If the value of $x_s$ remains unchanged or almost does not change
between the time layers~$n - 1$ and~$n + 1$, the numerical calculation of expression~(\ref{hardTerm})
causes practical difficulties.
Linearizing this expression and representing it in the form of an iterative process only complicates the situation,
so here we consider~$\Gamma^h$ without regard to the iterative process.
Notice that if
\[
|\hat{x}_s/ \check{x}_s| = 1 + \epsilon,
\quad
|\epsilon| \ll 1,
\]
then the following expansion can be considered instead of~(\ref{hardTerm})
\begin{equation}
\displaystyle
\frac{1}{\hat{x}_s - \check{x}_s} \, \ln \frac{\hat{x}_s}{\check{x}_s}
= \frac{\ln (1 + (\hat{x}_s/\check{x}_s-1))}{\check{x}_s(\hat{x}_s/\check{x}_s - 1)}
\sim
\frac{1}{\check{x}_s} \sum_{\kappa=0}^{\infty} \frac{1}{\kappa+1}\left(1 - \frac{\hat{x}_s}{\check{x}_s}\right)^\kappa.
\end{equation}
In the regions of small change in $x_s$, the latter expansion is used, and the first eight terms of the expansion are taken.
Outside such regions, calculations can still be performed using equation~(\ref{hardTerm}).

\subsubsection{Dam break over a parabolic bottom}

The dam-break problem is considered over a parabolic bottom
\begin{equation}
H(x) = d_1\left[ \left(\frac{2}{L}\right)^2 \left(x - \frac{L}{2}\right)^2 - 1\right]
\end{equation}
where $L$ is the length of the river segment and $d_1$ is the height of the parabolic bottom at the point~$x = L/2$.
According to~\cite{bk:DorKapJMPSW2021}, the following approximation for~$H^\prime(x)$ is chosen
\begin{equation}
(H^\prime)^h = \frac{2 (\cosh(\sqrt{\beta} \tau) - 1)}{\tau^2} \left(x -\frac{L}{2}\right),
\end{equation}
where $\beta = 8 d_1/L^2$.

\begin{figure}[H]
  \centering
  \includegraphics[width=0.4\linewidth]{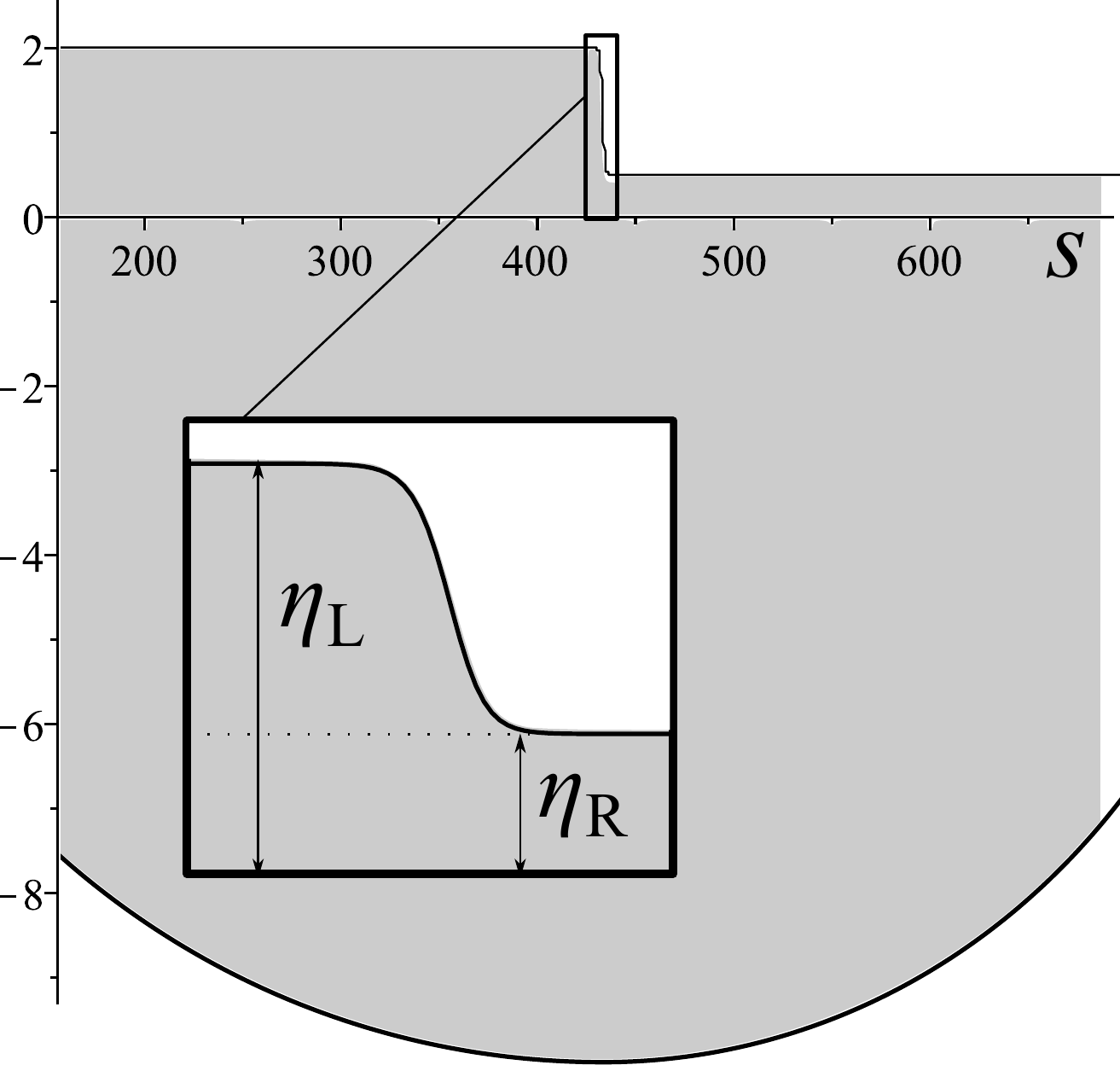}
  \caption{Initial water profile over the parabolic bottom (in Lagrangian coordinates). The smoothen initial data segment is shown on a larger scale.}
  \label{pic:ics1}
\end{figure}

In order to provide smoother initial data, the initial free surface profile is described by the function
\begin{equation}
\eta(\xi) = \eta_L - \frac{\eta_L - \eta_R}{1 + \exp\left(\sigma(\xi  -L/2)\right)},
\end{equation}
where $\sigma = 20$ is the curve steepness coefficient, and the constants~$\eta_L=2$ and~$\eta_R=0.5$
are given in~Figure~\ref{pic:ics1}. We also put $d_1=10$ and $L=100$.
By means of~(\ref{massLagrInt}) one states that the total mass of the fluid has value~$s \sim 791.7$.
Here and further we choose~$h=0.1$ and~$\tau=0.01$.

\begin{figure}[H]
  \centering
  \includegraphics[width=0.4\linewidth]{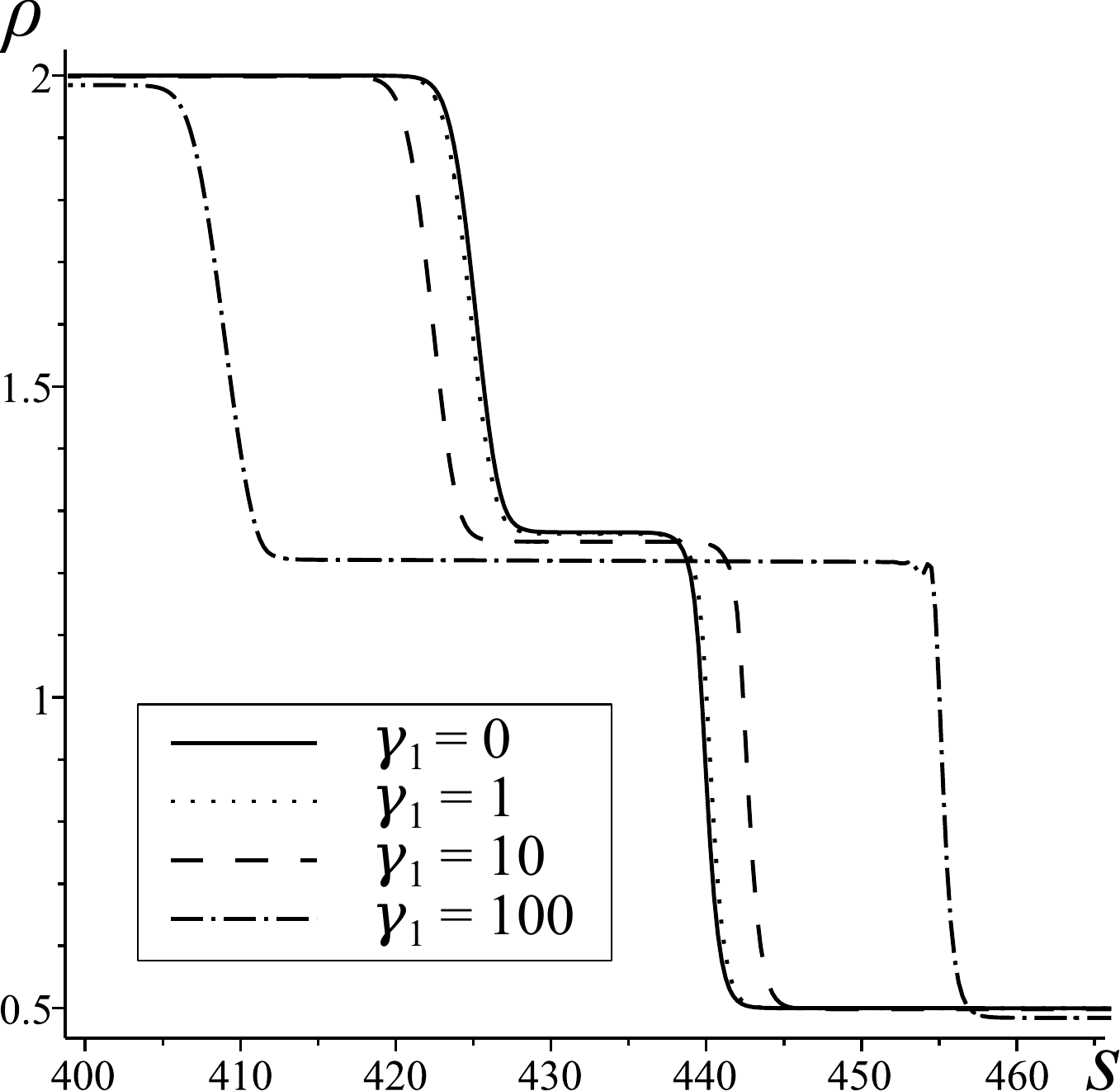}
  \caption{The approximate change in the solution profile at $t=0.2$ depending on~$\gamma_1$ calculated using scheme~(\ref{naive}). Case $\gamma_1=0$ corresponds to the one-dimensional shallow water equations.}
  \label{pic:gamma1cmp}
\end{figure}

To investigate the qualitative effect of the value of~$\gamma_1$ on the solution,
calculations have been performed at the moment~$t = 0.2$ for different values of~$\gamma_1$.
Figure~\ref{pic:gamma1cmp} shows that an increase in~$\gamma_1$ leads to
an approximately linear increase in the fluid velocity. Based on Figure~\ref{pic:gamma1cmp}, further on we put~$\gamma_1=10$.

\begin{figure}[H]
  \centering
  \includegraphics[width=0.6\linewidth]{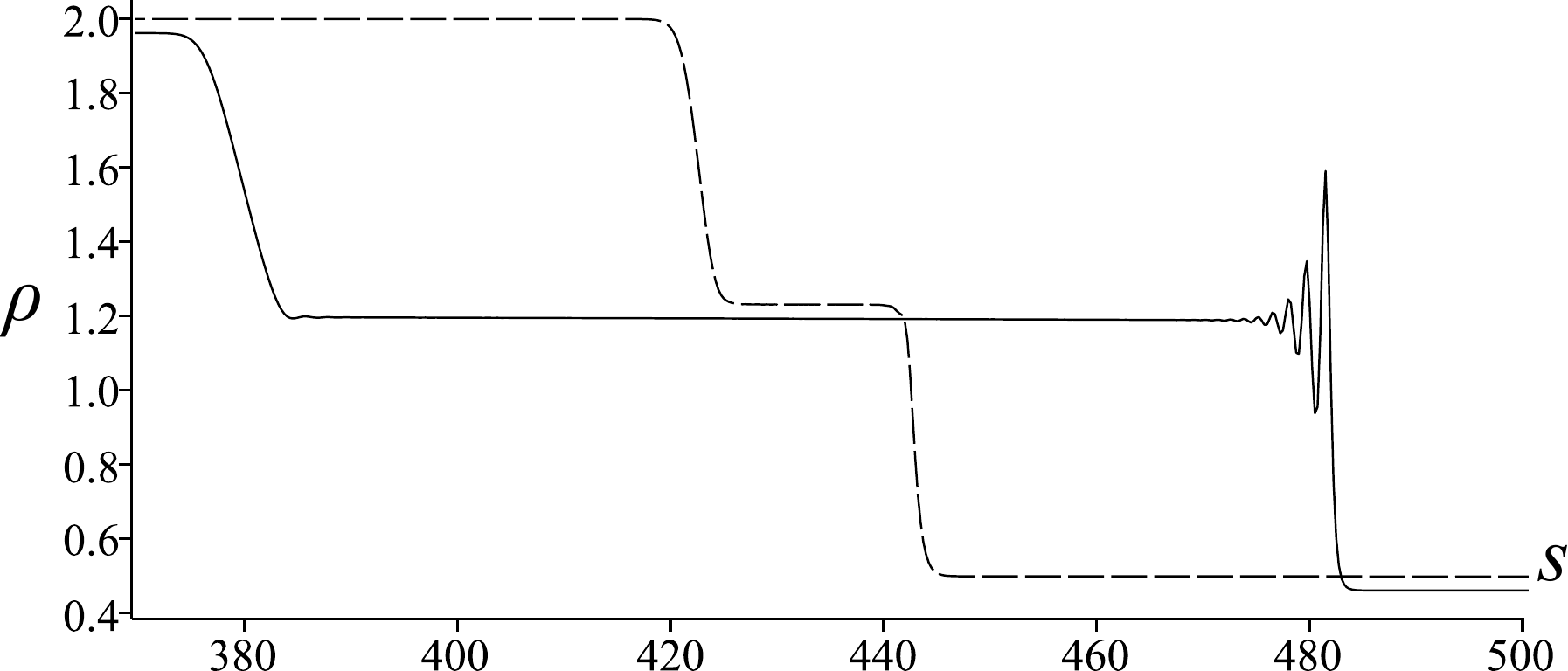}
  \caption{Solutions for~$\gamma_1=10$, $t_1=0.2$ (dash line) and $t_2=1$ (solid line).}
  \label{pic:sol1}
\end{figure}

\begin{figure}[H]
  \centering
  \includegraphics[width=0.6\linewidth]{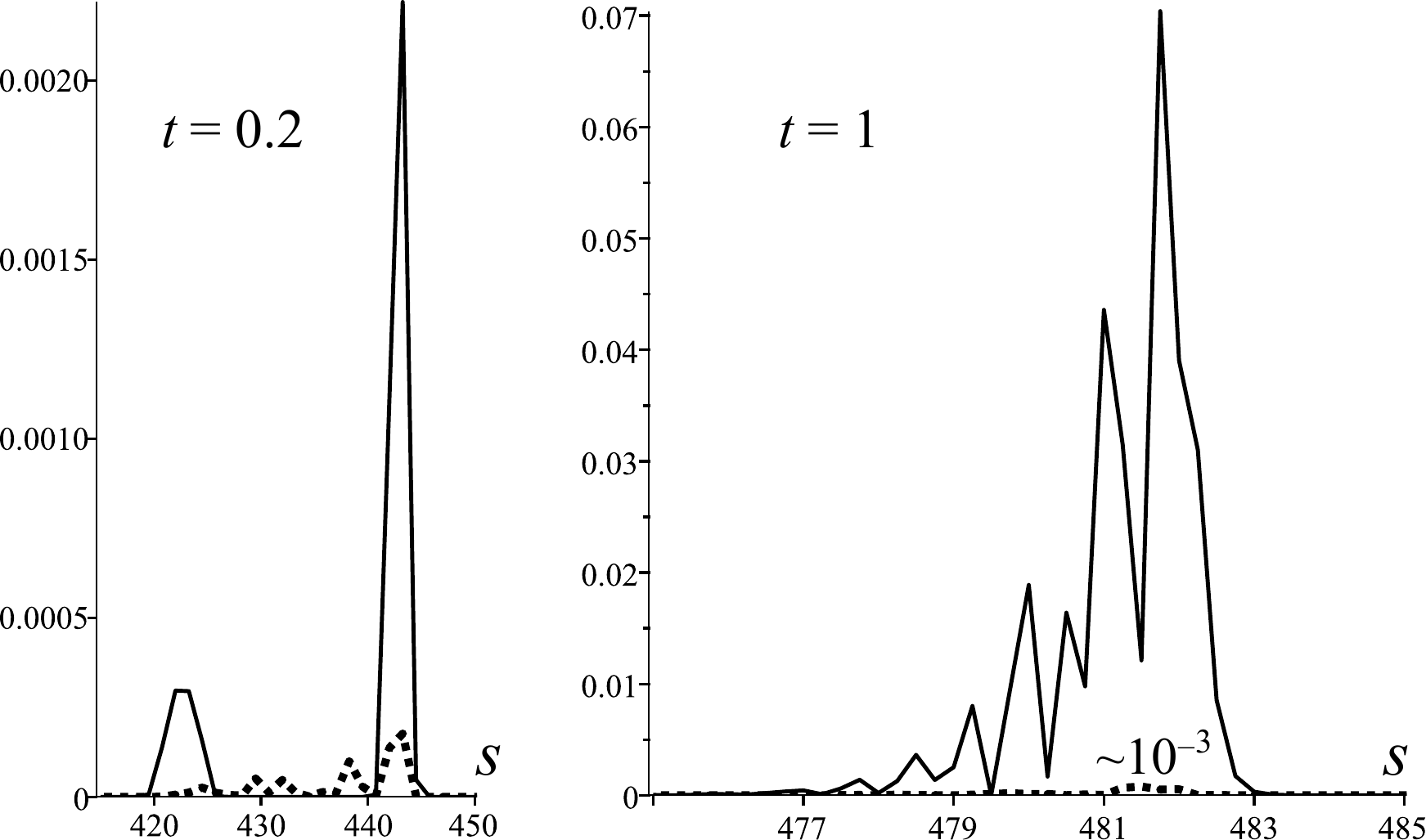}
  \caption{Energy conservation law errors comparison
  for the conservative scheme (dot line) and for ``naive'' scheme (solid line)
  at time $t_1=0.2$ and $t_2=1$.}
  \label{pic:err1}
\end{figure}

The solutions of the problem for $t_1=0.2$ and $t_2=1$ are given in Figure~(\ref{pic:sol1}).
Here and further no artificial viscosity is used since this allows better control of the energy preservation on solutions.
The profiles of solutions obtained by schemes~(\ref{naive}) and~(\ref{scheme}) on the given scale practically do not differ, therefore, solutions obtained by scheme~(\ref{scheme}) are presented throughout the text.
In contrast to solution profiles, the quality of energy conservation on solutions varies considerably for the schemes.
In Figure~\ref{pic:err1}, the energy preservation on solutions is given using the conservation law~(\ref{SchemeCLenergyHParabolic}) and the estimation~(\ref{deltaEpsExpr}).
Scheme~(\ref{scheme}) conserves energy much better than ``naive'' scheme (\ref{naive}).

\begin{figure}[H]
  \centering
  \includegraphics[width=0.4\linewidth]{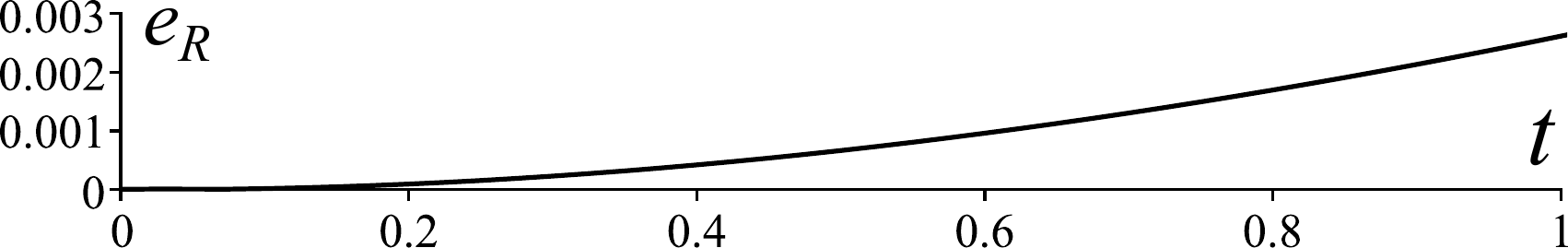}
  \caption{Total energy relative error for the schemes on time interval $0 \leqslant t \leqslant 1$.}
  \label{pic:err1ham}
\end{figure}

The estimates of the total energy conservation by formula~(\ref{ErrRelativeFormula}) for the two schemes practically do not differ as it shown in Figure~\ref{pic:err1ham}. The total energy conservation pattern will be significantly different for the next problem.

\subsubsection{Collapse of a fluid column above an inclined bottom}

\begin{figure}[H]
  \centering

  \begin{subfigure}[b]{0.415\linewidth}
    \includegraphics[width=\linewidth]{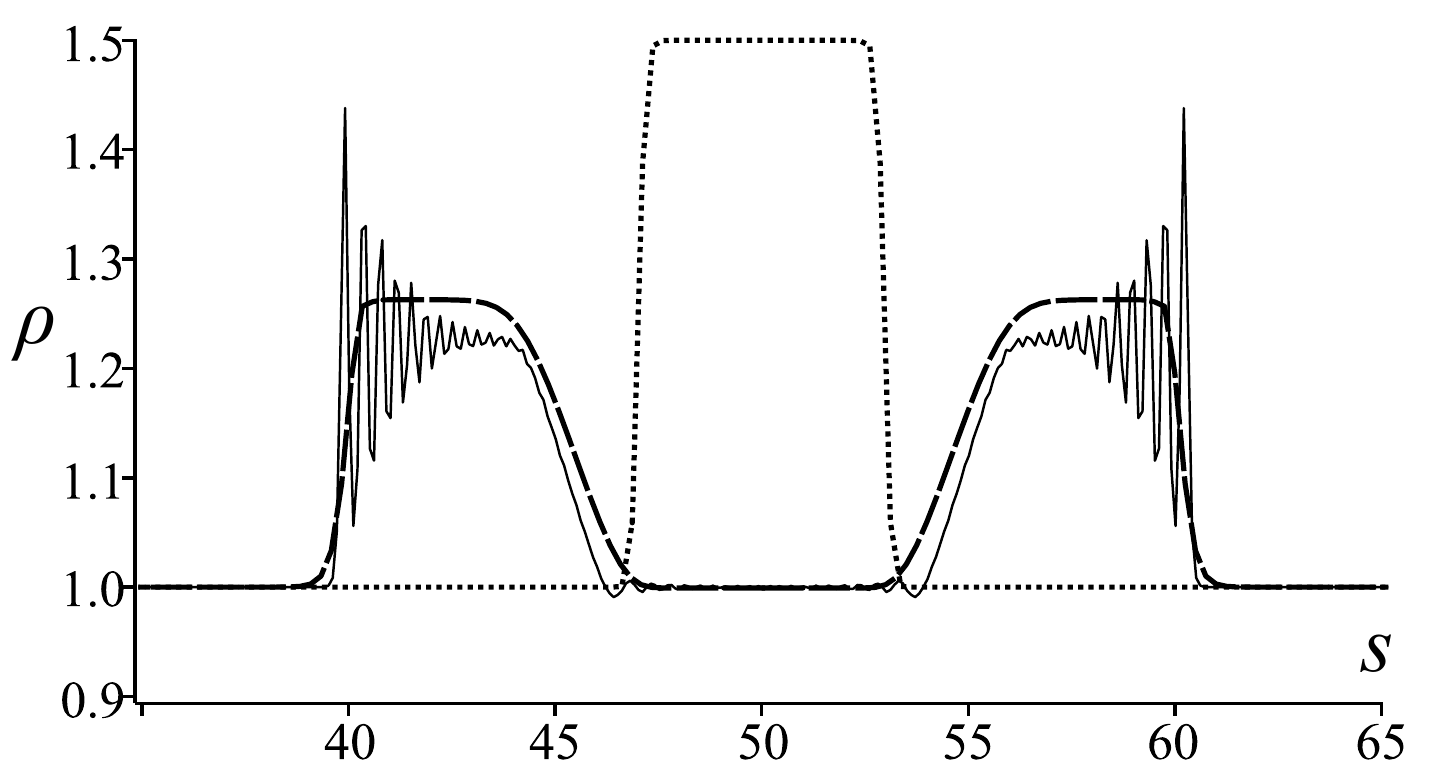}
  \end{subfigure}

  \vfill

  \begin{subfigure}[b]{0.4\linewidth}
    \includegraphics[width=\linewidth]{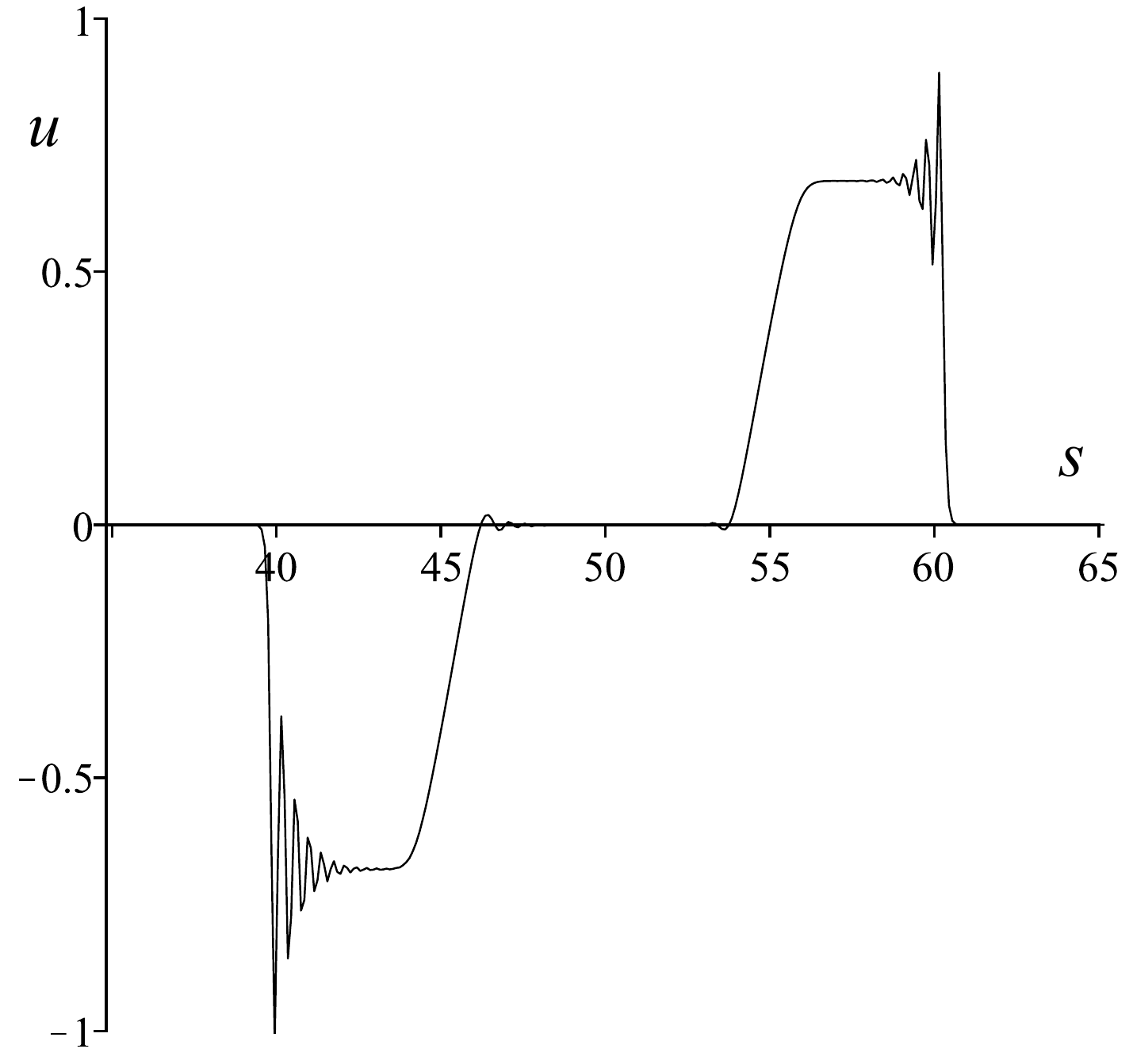}
  \end{subfigure}

  \caption{Height of the fluid (top) and velocity (bottom) at~$t=2$ are depicted with solid lines.
  The dot line is used for initial data profile, and the dash line is used for the solution with artificial viscosity.}
  \label{pic:sol2}
\end{figure}

A column of liquid above an inclined bottom, which collapses due to gravity, is considered.
By means of transformation~(\ref{InclToFlatBottomTrDelta}), the problem is reduced to the case of a horizontal bottom.
The initial data is depicted in Figure~\ref{pic:sol2} by dot line. The initial profile is smoothed in the same way as in the previous problem. It is described by the function
\begin{equation}
\eta(\xi) = \eta_L
    - \frac{\eta_L - \eta_R}{1 + \exp\left(\sigma(\xi - L/2 + dL)\right)}
    + \frac{\eta_L - \eta_R}{1 + \exp\left(\sigma(\xi - L/2 - dL)\right)},
\end{equation}
where $dL=2$, and the remaining parameters have the same meaning and values as in the previous section.

\begin{figure}[H]
  \centering
  \includegraphics[width=0.4\linewidth]{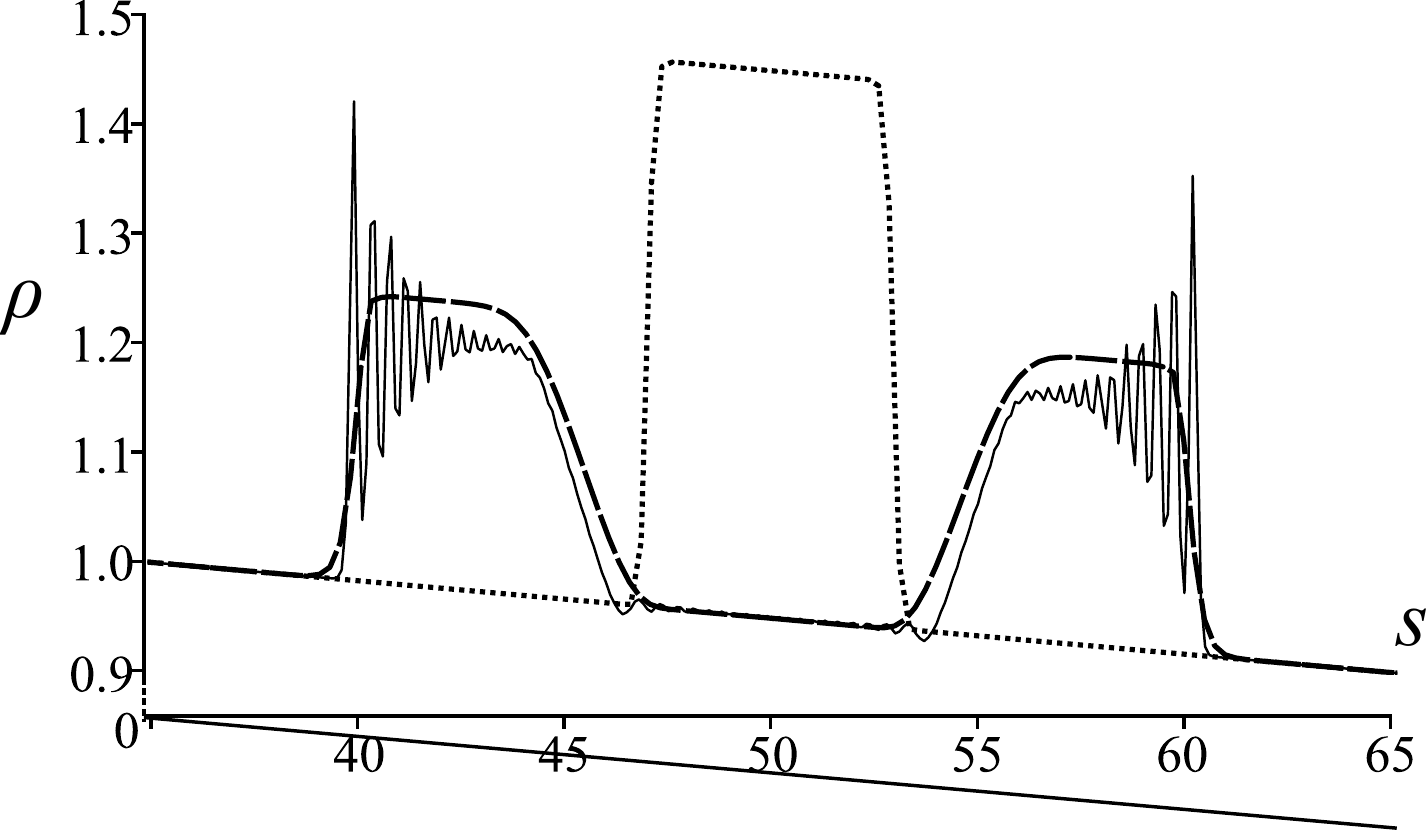}
  \caption{The solution (see Figure~\ref{pic:sol2}) for an inclined bottom obtained
  by means of transformation~(\ref{InclToFlatBottomTrDelta}).
  The solution was also shifted along the vertical axis for clarity.}
  \label{pic:sol2skew}
\end{figure}

\begin{figure}[H]
  \centering
  \includegraphics[width=0.5\linewidth]{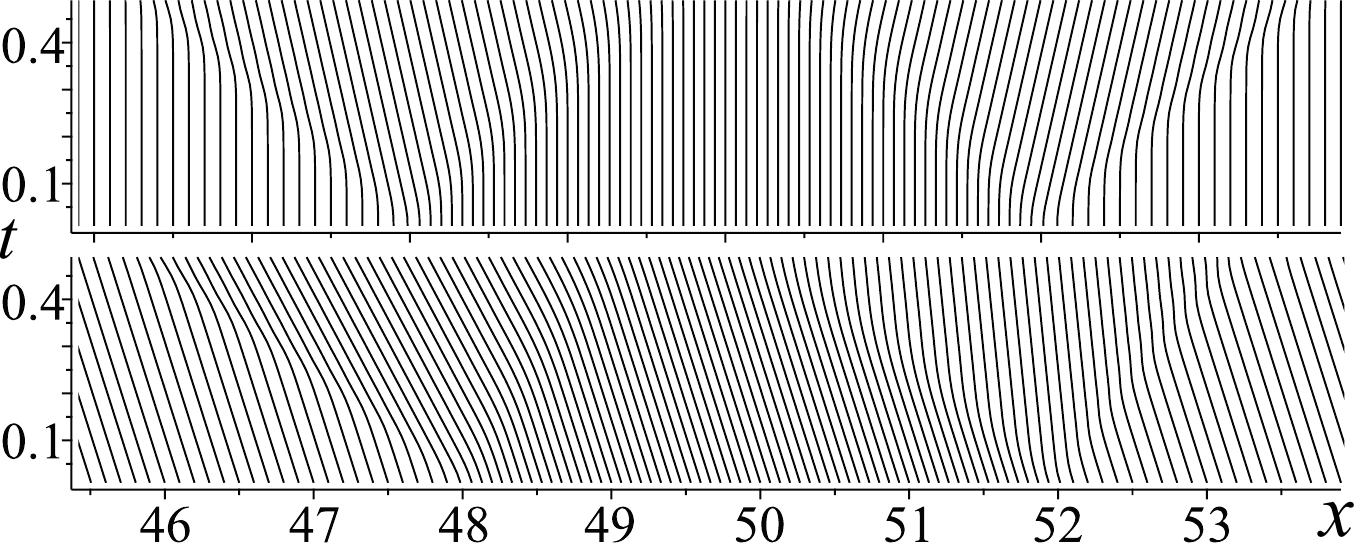}
  \caption{Top to bottom: trajectories for~$u_0=0$ and~$u_0=-1$ on the time interval~$0 \leqslant t \leqslant 0.5$.}
  \label{pic:sol2x}
\end{figure}

In Figure~\ref{pic:sol2}, the solid line shows the solution at time $t = 2$.
The dashed line shows the solution profile calculated using artificial viscosity.
The solution for an inclined bottom is given in Figure~\ref{pic:sol2skew}.
This solution is very similar to the solution obtained in Eulerian coordinates in~\cite{bk:KarPetSla09}.
In~\cite{bk:KarPetSla09}, the case is also considered when the liquid rises at a certain velocity~$u_0<0$ up the inclined plane.
In contrast to Eulerian coordinates, in Lagrangian coordinates the profiles of solutions for this case do not differ from the case with zero velocity, since the mass distribution remains the same.
In these cases, only the trajectories of the particles differ as shown in Figure~\ref{pic:sol2x}.

\begin{figure}[H]
  \centering
  \includegraphics[width=0.4\linewidth]{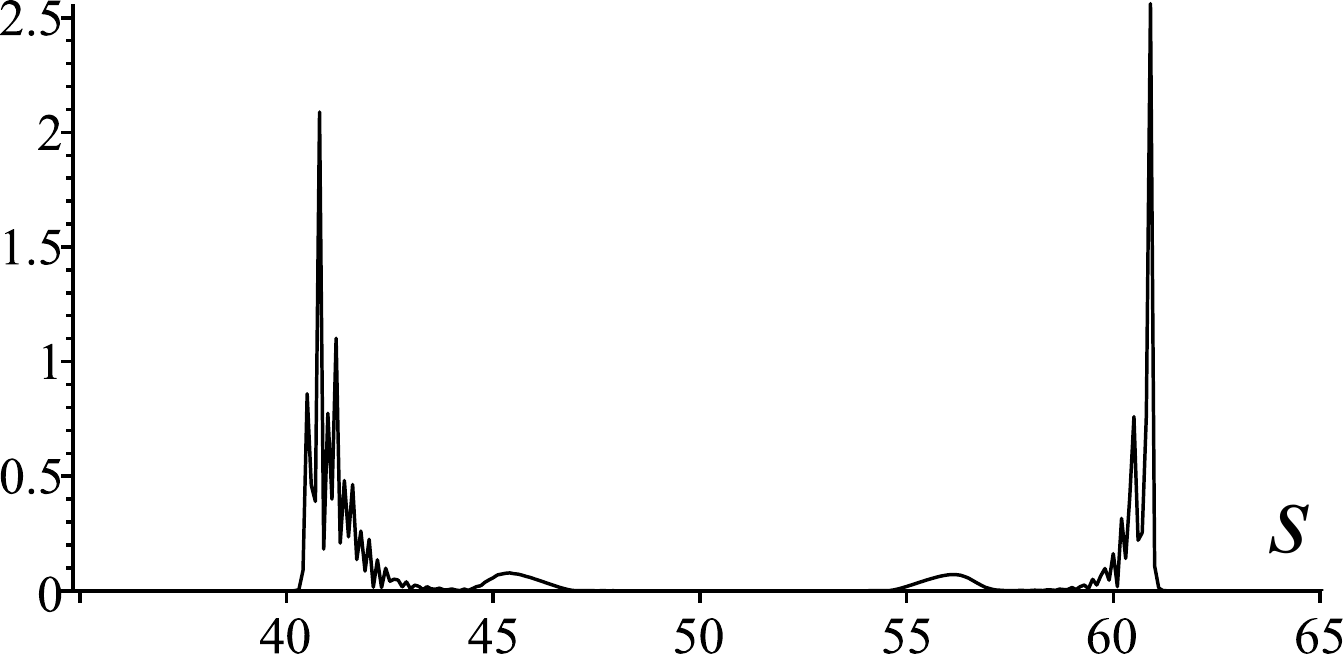}
  \hfil
  \includegraphics[width=0.4\linewidth]{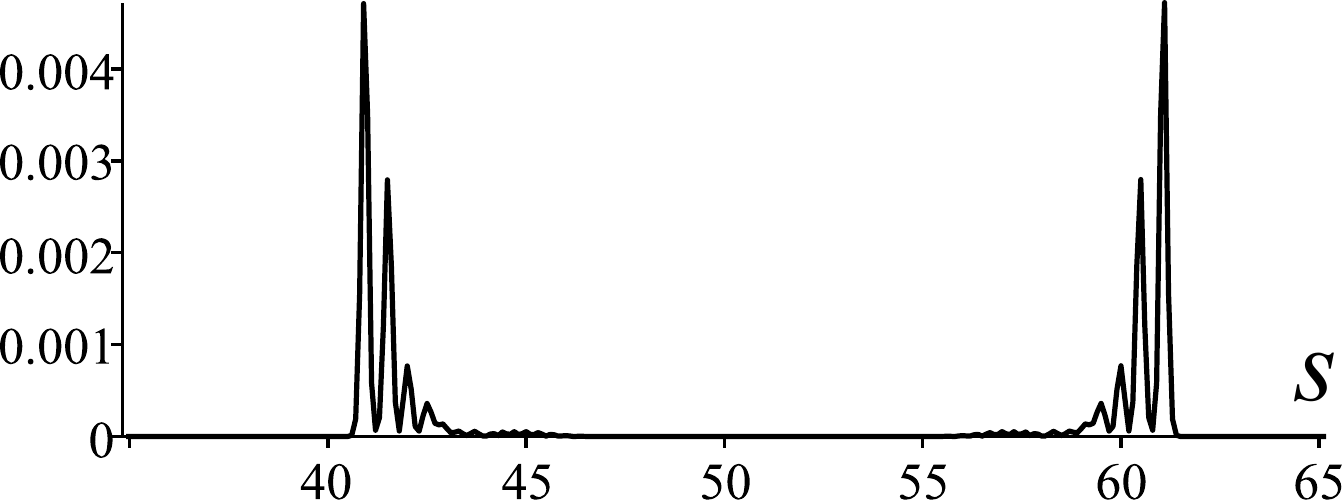}
  \caption{Left to right: energy conservation law errors for ``naive'' scheme~(\ref{naive}) and conservative scheme~(\ref{scheme}) at~$t=2$.}
  \label{pic:err2}
\end{figure}

Figure~\ref{pic:err2} is similar to Figure~\ref{pic:err1} from the previous section. It again demonstrates a significant difference in energy preservation between schemes~(\ref{scheme}) and~(\ref{naive}).

\begin{figure}[H]
  \centering
  \includegraphics[width=0.7\linewidth]{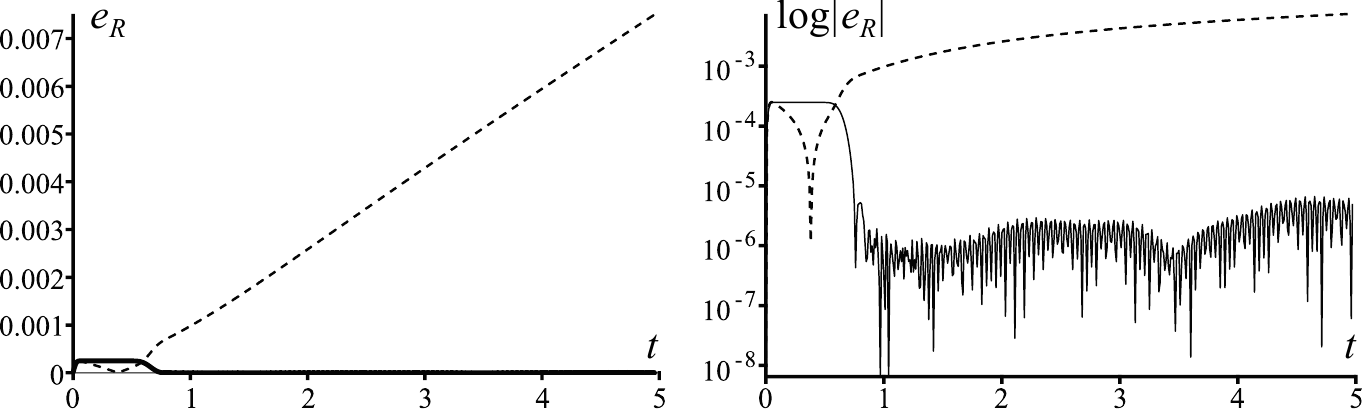}
  \caption{Total energy relative error (left) and its logarithm (right) on the time interval~$0 \leqslant t \leqslant 5$ for ``naive`` scheme (dash line) and the conservative scheme (solid line).}
  \label{pic:err2ham}
\end{figure}

Figure~\ref{pic:err2ham} shows that the increase in the error~(\ref{ErrRelativeFormula}) over time for ``naive'' scheme is quite significant in contrast to the slowly growing error for the conservative scheme.

\section{Conclusion}

Symmetries and conservation laws of the one-dimensional modified shallow water equations
 in Lagrangian coordinates for various bottom topographies are considered.
These results are based on group classification
\cite{bk:KaptsovMeleshko_1D_classf[2018]}. Variational formulation
of equations in Lagrangian coordinates allows one to obtain their
conservation laws by means of the Noether theorem. The corresponding
conservation laws in mass Lagrangian coordinates and Eulerian
coordinates derived from conservation laws in Lagrangian coordinates
are also given.

On the basis of invariant finite-difference schemes for the shallow
water equations recently obtained
in~\cite{dorodnitsyn2019shallow,bk:DorKapJMPSW2021}, invariant
schemes for the modified shallow water equations for various bottom
topographies are constructed. These schemes possess
finite-difference analogues of the local conservation laws of mass and
energy for arbitrary shape of bottom, as well as additional
conservation laws that appear for special cases of bottom
topography. All the schemes are constructed on uniform orthogonal
meshes which are invariant with respect to all symmetries inherited
from the differential model.

To construct conservative schemes on such meshes, it is often
convenient to use the finite-difference analogue of the  direct
method. The direct method is well suited for schemes that can be
written in terms of rational expressions. This approach successfully
used to construct conservative schemes for the standard shallow water
equations, but for the schemes for the modified shallow water
equations it doesn't work by straightforward application. Indeed, in
the case of the modified shallow water equations the conservation
laws of energy and momentum include logarithmic terms which
significantly complicates the problem. For the problem under
consideration, the authors managed to find a special ansatz which
solves the problem of logarithmic term. This example shows the
importance of such methods as the finite-difference analogue of the
Noether
 theorem~\cite{Dor_3,bk:Dorod_NoetherTh[2001],bk:DorodKozlovWint[2004]}
and the Lagrange identity and adjoint equation
method~\cite{bk:DorodKozlovWintKaptsov[2014],bk:DorodKozlovWintKaptsov[2015]}.

The numerical implementation of the constructed finite-difference
schemes is carried out for the examples of a dam break over a parabolic bottom
and a collapsing liquid column over an inclined bottom. All calculations
are performed in Lagrangian coordinates. The constructed schemes are
compared with a naive invariant scheme constructed without invariant
Lagrangian consideration. It is shown that the specially constructed
conservative schemes preserve energy much better than the naive
approximation. This emphasizes the importance of the criteria of
invariance and conservativeness in the construction of schemes.

\section*{Acknowledgements}

%\todo{Anyway, for arXiv it should be OK. Discuss with S.V. before submitting to a journal!}

The research  was supported by Russian Science Foundation Grant No
18-11-00238 ``Hy\-dro\-dy\-nam\-ics-type equations: symmetries, conservation
laws, invariant difference schemes''. S.M. and E.K. thanks Suranaree
University of Technology (SUT) and  Thailand Science Research and Innovation (TSRI).

%%%%%%%%%%%%%%%%%%%%%%%%%%%%%%%%%%%%%%%%%%%%%%%%%%%%%%%%%%%%%%%%%%%%%%%%%%%%


\begin{thebibliography}{99}

\bibitem{bk:BonnetoBarthelemy_et[2011]}
P.~Bonneton, E.~Barth\'elemy, F.~Chazel, R.~Cienfuegos, D.~Lannes, F.~Marche,
  and M.~Tissier.
\newblock Recent advances in {S}erre-{G}reen-{N}aghdi modelling for wave
  transformation, breaking and runup processes.
\newblock {\em Euro. J. Mech. B/Fluids}, 30:589--597, 2011.

\bibitem{bk:KhakimzyanovDutykhFedotovaMitsotakis}
G.~S. Khakimzyanov, D.~Dutykh, Z.~I. Fedotova, and D.~E. Mitsotakis.
\newblock Dispersive shallow water wave modelling. {P}art {I}: {M}odel
  derivation on a globally flat space.
\newblock {\em Commun. Comput. Phys.}, 23(1):1--29, 2018.

\bibitem{bk:Whitham[1974]}
G.~B. Whitham.
\newblock {\em Linear and Nonlinear Waves}.
\newblock Wiley, New York, 1974.

\bibitem{bk:Ovsyannikov[2003]}
L.~V. Ovsiannikov.
\newblock {\em Lectures on the gas dynamics equations}.
\newblock Institute of Computer Studies, Moscow--Izhevsk, 2003.
\newblock in Russian.

\bibitem{bk:PetrosyanBook[2010]}
A.~S. Petrosyan.
\newblock {\em Additional chapters of heavy fluid hydrodynamics with a free
  boundary}.
\newblock Space Research Institute of the Russian Academy of Sciences, Moscow,
  2014.
\newblock in Russian.

\bibitem{bk:Vallis[2006]}
G.~K. Vallis.
\newblock {\em Atmospheric and Oceanic Fluid Dynamics: Fundamentals and
  Large-scale Circulation}.
\newblock Cambridge University Press, Cambridge, 2006.

\bibitem{bk:LeVeque}
R.~J. LeVeque.
\newblock {\em Finite Volume Methods For Hyperbolic Problems. Cambridge Texts
  in Applied Mathematics}.
\newblock Cambridge University Press, Cambridge, 2002.

\bibitem{bk:Pedlosky}
J.~Pedlosky.
\newblock {\em Geophysical Fluid Dynamics}.
\newblock Springer, New York, 1987.

\bibitem{bk:KarelskyPetrosyan2006}
K.~V. Karelsky and A.~S. Petrosyan.
\newblock Particular solutions and {R}iemann problem for modified shallow water
  equations.
\newblock {\em Fluid Dynamics Research}, 38(5):339--358, 2006.

\bibitem{bk:KarelskyPetrosyan2008}
K.~V. Karelsky and A.~S. Petrosyan.
\newblock Modified shallow water equations. {S}imple waves and {R}iemann
  problem.
\newblock In I.S.~Mamaev A.V.~Borisov, V.V.~Kozlov and M.A. Sokolovskiy,
  editors, {\em IUTAM Symposium on Hamiltonian Dynamics, Vortex Structures,
  Turbulence}, pages 383--392. Springer, Dordrecht, 2008.
\newblock Proceedings of the IUTAM Symposium held in Moscow, 25--30, August,
  2006.

\bibitem{bk:Ovsyannikov[1982]}
L.~V. Ovsiannikov.
\newblock {\em Group Analysis of Differential Equations}.
\newblock Academic, New York, 1982.

\bibitem{bk:Olver[1986]}
P.~J. Olver.
\newblock {\em Applications of {L}ie Groups to Differential Equations}.
\newblock Springer-Verlag, New York, 1986.

\bibitem{bk:Ibragimov1985}
N.~H. Ibragimov.
\newblock {\em Transformation Groups Applied to Mathematical Physics}.
\newblock Reidel, Boston, 1985.

\bibitem{bk:Bluman1989}
G.W. Bluman and S.~Kumei.
\newblock {\em Symmetries and Differential Equations}.
\newblock Applied Mathematical Sciences. Springer New York, 2013.

\bibitem{bk:HandbookLie_v1}
N.~H. Ibragimov, editor.
\newblock {\em {CRC} Handbook of {L}ie Group Analysis of Differential
  Equations}, volume~1.
\newblock CRC Press, Boca Raton, 1994.
\newblock N. H. Ibragimov ({\it ed.}).

\bibitem{bk:Gaeta1994}
G.~Gaeta.
\newblock {\em Nonlinear Symmetries and Nonlinear Equations}.
\newblock Kluwer, Dordrecht, 1994.

\bibitem{bk:HandbookLie_v2}
N.~H. Ibragimov, editor.
\newblock {\em {CRC} Handbook of {L}ie Group Analysis of Differential
  Equations}, volume~2.
\newblock CRC Press, Boca Raton, 1995.

\bibitem{bk:LeviNicciRogersWint[1989]}
D.~Levi, M.C. Nicci, C.~Rogers, and P.~Winternitz.
\newblock Group theoretical analysis of a rotating shallow liquid in a rigid
  container.
\newblock {\em Journal of Physics A:~Mathematical and General}, 22:4743--4767,
  1989.

\bibitem{bk:ClarksonBila[2006]}
N.~Bila, E.~Mansfield, and P.~Clarkson.
\newblock Symmetry group analysis of the shallow water and semi-geostrophic
  equations.
\newblock {\em The Quarterly Journal of Mechanics and Applied Mathematics}, 59,
  02 2006.

\bibitem{bk:Andronikos2019}
A.~Paliathanasis.
\newblock Lie symmetries and similarity solutions for rotating shallow water.
\newblock {\em {Z}eitschrift für {N}aturforschung {A}}, 06 2019.

\bibitem{bk:AksenovDruzkov_classif[2019]}
A.~V. Aksenov and K.~P. Druzhkov.
\newblock Conservation laws of the equation of one-dimensional shallow water
  over uneven bottom in {L}agrange's variables.
\newblock {\em International Journal of Non-Linear Mechanics}, 119:103348,
  2020.

\bibitem{bk:KaptsovMeleshko_1D_classf[2018]}
E.~I. Kaptsov and S.~V. Meleshko.
\newblock Analysis of the one-dimensional {E}uler--{L}agrange equation of
  continuum mechanics with a {L}agrangian of a special form.
\newblock {\em Applied Mathematical Modelling}, 77:1497 -- 1511, 2020.

\bibitem{bk:SiriwatKaewmaneeMeleshko2016}
P.~Siriwat, C.~Kaewmanee, and S.~V. Meleshko.
\newblock Symmetries of the hyperbolic shallow water equations and the
  {G}reen-{N}aghdi model in {L}agrangian coordinates.
\newblock {\em International Journal of Non-Linear Mechanics}, 86:185--195,
  2016.

\bibitem{bk:Bernetti[2008]}
R.~Bernetti, V.~A. Titarev, and E.~F. Toro.
\newblock Exact solution of the riemann problem for the shallow water equations
  with discontinuous bottom geometry.
\newblock {\em Journal of Computational Physics}, 227(6):3212 -- 3243, 2008.

\bibitem{bk:HanHantke[2012]}
E.~E. Han, M.~Hantke, and G.~Warnecke.
\newblock Exact riemann solutions to compressible {E}uler equations in ducts
  with discontinuous cross-section.
\newblock {\em Journal of Hyperbolic Differential Equations}, 09(03):403--449,
  2012.

\bibitem{Maeda1}
S.~Maeda.
\newblock Extension of discrete {Noether} theorem.
\newblock {\em Math. Japonica}, 26(1):85--90, 1985.

\bibitem{Maeda2}
S.~Maeda.
\newblock The similarity method for difference equations.
\newblock {\em J. Inst. Math. Appl.}, 38:129--134, 1987.

\bibitem{Dor_1}
V.~A. Dorodnitsyn.
\newblock Transformation groups in net spaces.
\newblock {\em Journal of Soviet Mathematics}, 55(1):1490--1517, Jun 1991.

\bibitem{Dor_2}
V.~A. Dorodnitsyn.
\newblock Finite difference models entirely inheriting symmetry of original
  differential equations.
\newblock {\em International Journal of Modern Physics C}, 5, 08 1994.

\bibitem{Dor_3}
V.~A. Dorodnitsyn.
\newblock The finite-difference analogy of {Noether}'s theorem.
\newblock {\em Phys. Dokl.}, 38(02):66--68, 1993.

\bibitem{bk:DorodKozlovWint[2004]}
V.~A. Dorodnitsyn, R.~V. Kozlov, and P.~Winternitz.
\newblock Continuous symmetries of {L}agrangians and exact solutions of
  discrete equations.
\newblock {\em Journal of Mathematical Physics}, 45(1):336--359, 2004.

\bibitem{KimOlv04}
Kim P. and P.~J. Olver.
\newblock Geometric integration via multi-space.
\newblock {\em Regul. Chaotic Dyn.}, 9:213--226, 2004.

\bibitem{KimOlv07}
M.~Welk, P.~Kim, and P.J. Olver.
\newblock Numerical invariantization for morphological {PDE} schemes.
\newblock In F.~Sgallari, A.~Murli, and N.~Paragios, editors, {\em Scale Space
  and Variational Methods in Computer Vision}, pages 508--519, Berlin,
  Heidelberg, 2007. Springer Berlin Heidelberg.

\bibitem{[LW-2]}
D.~Levi and P.~Winternitz.
\newblock Continuous symmetries of difference equations.
\newblock {\em Journal of Physics A: Mathematical and General}, 39(2):R1--R63,
  12 2005.

\bibitem{bk:DorodKozlovWinternitz[2000]}
V.~A. Dorodnitsyn, R.~V. Kozlov, and P.~Winternitz.
\newblock Lie group classification of second-order ordinary difference
  equations.
\newblock {\em Journal of Mathematical Physics}, 41(1):480--504, 2000.

\bibitem{bk:Dorod_NoetherTh[2001]}
V.~A. Dorodnitsyn.
\newblock {Noether}-type theorems for difference equations.
\newblock {\em Applied Numerical Mathematics}, 39(3):307 -- 321, 2001.
\newblock Themes in Geometric Integration.

\bibitem{Quisp}
G.R.W. Quispel and R.~Sahadevan.
\newblock Lie symmetries and the integration of difference equations.
\newblock {\em Physics Letters A}, 184(1):64 -- 70, 1993.

\bibitem{[LW-3]}
P.~Winternitz.
\newblock Symmetry preserving discretization of differential equations and
  {L}ie point symmetries of differential-difference equations.
\newblock pages 292--341, 2011.

\bibitem{bk:Dorodnitsyn[2011]}
V.~A. Dorodnitsyn.
\newblock {\em Applications of Lie Groups to Difference Equations}.
\newblock CRC Press, Boca Raton, 2011.

\bibitem{Vinet}
Roberto Floreanini and Luc Vinet.
\newblock Lie symmetries of finite‐difference equations.
\newblock {\em Journal of Mathematical Physics}, 36(12):7024--7042, 1995.

\bibitem{bk:Hydon_book[2014]}
P.~E. Hydon.
\newblock {\em Difference Equations by Differential Equation Methods}.
\newblock Cambridge Monographs on Applied and Computational Mathematics.
  Cambridge University Press, 2014.

\bibitem{bk:DorodKozlovWintKaptsov[2015]}
V.~A. Dorodnitsyn, E.~I. Kaptsov, R.~V. Kozlov, and P.~Winternitz.
\newblock The adjoint equation method for constructing first integrals of
  difference equations.
\newblock {\em Journal of Physics A: Mathematical and Theoretical},
  48(5):055202, 01 2015.

\bibitem{[Dheat]}
V.~A. Dorodnitsyn and R.~V. Kozlov.
\newblock A heat transfer with a source: the complete set of invariant
  difference schemes.
\newblock {\em Journal of Nonlinear Mathematical Physics}, 10, 10 2003.

\bibitem{bk:DorodKaptsov[2013]}
V.~A. Dorodnitsyn and E.~I. Kaptsov.
\newblock Discretization of second-order ordinary differential equations with
  symmetries.
\newblock {\em Computational Mathematics and Mathematical Physics},
  53(8):1153--1178, 2013.

\bibitem{bk:DorodKaptsov_Ermakov[2016]}
V.~A. Dorodnitsyn and E.~I. Kaptsov.
\newblock Invariant difference schemes for the {E}rmakov system.
\newblock {\em Differential Equations}, 52:926--941, 01 2016.

\bibitem{bk:Dorod_Hamilt[2011]}
V.~A. Dorodnitsyn and R.~V Kozlov.
\newblock {L}agrangian and {H}amiltonian formalism for discrete equations:
  Symmetries and first integrals.
\newblock In D.~Levi, P.~Olver, Z.~Thomova, and P.~Winternitz, editors, {\em
  Symmetries and Integrability of Difference Equations}, London Mathematical
  Society Lecture Note Series, page 7--49. Cambridge University Press, 2011.

\bibitem{bk:Dorod_Hamilt[2010]}
V.~A. Dorodnitsyn and R.~V. Kozlov.
\newblock Invariance and first integrals of continuous and discrete
  {H}amiltonian equations.
\newblock {\em Journal of Engineering Mathematics}, 66(1):253--270, Mar 2010.

\bibitem{bk:BlumanAnco2002}
G.~W. Bluman and S.~C. Anco.
\newblock {\em Symmetry and Integration Methods for Differential Equations}.
\newblock Springer, New York, 2002.

\bibitem{bk:DorodKozlovWintKaptsov[2014]}
P.~Winternitz, V.~A. Dorodnitsyn, E.~I. Kaptsov, and R.~V. Kozlov.
\newblock First integrals of difference equations which do not possess a
  variational formulation.
\newblock {\em Doklady Mathematics}, 89(1):106--109, 01 2014.

\bibitem{bk:Bluman1997}
S.~C. Anco and G.~W. Bluman.
\newblock Direct construction of conservation laws from field equations.
\newblock {\em Physical Review Letters}, 78:2869--2873, 04 1997.

\bibitem{bk:BlumanCheviakovAnco}
G.~W. Bluman, A.~F. Cheviakov, and S.~C. Anco.
\newblock {\em Applications of Symmetry Methods to Partial Differential
  Equations}.
\newblock Springer, New York, 2010.
\newblock Applied Mathematical Sciences, Vol.168.

\bibitem{dorodnitsyn2019shallow}
V.~A. Dorodnitsyn and E.~I. Kaptsov.
\newblock Shallow water equations in {L}agrangian coordinates: {S}ymmetries,
  conservation laws and its preservation in difference models.
\newblock {\em Communications in Nonlinear Science and Numerical Simulation},
  89:105343, 2020.

\bibitem{bk:DorKapJMPSW2021}
V.~A. Dorodnitsyn and E.~I. Kaptsov.
\newblock Discrete shallow water equations preserving symmetries and
  conservation laws.
\newblock {\em Journal of Mathematical Physics}, 62(8):083508, 2021.

\bibitem{bk:Noether1918}
E.~Noether.
\newblock Invariante variations problem.
\newblock {\em Konigliche Gesellschaft der Wissenschaften zu Gottingen,
  Nachrichten, Mathematisch-Physikalische Klasse Heft 2}, pages 235--257, 1918.
\newblock English translation:~Transport Theory and Statist. Phys., 1(3), 1971,
  183-207.

\bibitem{bk:SamarskyPopov_book[1992]}
A.~A. Samarskii and Y.~P. Popov.
\newblock {\em Difference methods for solving problems of gas dynamics}.
\newblock Nauka, Moscow, 1980.
\newblock In Russian.

\bibitem{bk:YanenkRojd[1968]}
B.~L. Rojdestvenskiy and N.~N. Yanenko.
\newblock {\em Systems of quasilinear equations and their applications to gas
  dynamics}.
\newblock Nauka, Moscow, 1968.
\newblock in Russian.

\bibitem{bk:DorodKapMelGN2020}
V.A. Dorodnitsyn, E.I. Kaptsov, and S.V. Meleshko.
\newblock Symmetries, conservation laws, invariant solutions and difference
  schemes of the one-dimensional {G}reen--{N}aghdi equations.
\newblock {\em Journal of Nonlinear Mathematical Physics}, 28:90--107, 2020.

\bibitem{bk:DorKapSW2020}
V.~A. Dorodnitsyn and E.~I. Kaptsov.
\newblock Invariant conservative difference schemes for shallow water equations
  in {E}ulerian and {L}agrangian coordinates.
\newblock {\em Communications in Nonlinear Science and Numerical Simulation},
  2020.
\newblock Submitted.

\bibitem{bk:DorKapMelSW2DarXiv}
V.~A. {Dorodnitsyn}, E.~I. {Kaptsov}, and S.~V. {Meleshko}.
\newblock Symmetries, conservation laws and difference schemes of the
  (1+2)-dimensional shallow water equations in {L}agrangian coordinates.
\newblock {\em arXiv e-prints}, page arXiv:2012.04410, 2020.

\bibitem{bk:ChevDorKap2020}
A.~F. Cheviakov, V.~A. Dorodnitsyn, and E.~I. Kaptsov.
\newblock Invariant conservation law-preserving discretizations of linear and
  nonlinear wave equations.
\newblock {\em Journal of Mathematical Physics}, 2020.
\newblock Submitted.

\bibitem{bk:Korobitsyn_scheme[1989]}
V.~A. Korobitsyn.
\newblock Thermodynamically matched difference schemes.
\newblock {\em U.S.S.R. Comput. Math. Math. Phys}, 29:71--79, 1989.

\bibitem{bk:KarPetSla09}
K.~V. Karelsky, A.~S. Petrosyan, and A.~G. Slavin.
\newblock Study of shallow water flows over an arbitrary bed profile in the
  presence of external force.
\newblock {\em Matem. Mod.}, 21:41--58, 2009.

\bibitem{bk:Samarskii2001theory}
A.~A. Samarskiy.
\newblock {\em The Theory of Difference Schemes}.
\newblock Monographs and textbooks in pure and applied mathematics. CRC Press,
  2001.

\end{thebibliography}
\end{document}